\documentclass[12pt]{amsart}
\usepackage{epsfig,apu,vmargin,amsmath}
\usepackage[T1]{fontenc}
\usepackage[english]{babel}
\usepackage{ae,aecompl}  

\setpapersize{A4}  
\setmarginsrb{30mm}{20mm}{25mm}{30mm}{12pt}{9mm}{0pt}{5mm}

\begin{document}
\title[simulation of multibody systems]{On the simulation of multibody systems with holonomic constraints}
\author{Jukka Tuomela}
\address{Department of Mathematics,
     University of Joensuu,
     PL 111, 80101 Joensuu, Finland}
\email{jukka.tuomela@joensuu.fi}
\author{Teijo Arponen}
\address{Institute of Mathematics,
     Helsinki University of Technology,
     PL 1100, 02015 TKK, Finland}
\email{teijo.arponen@hut.fi}
\author{Villesamuli Normi}
\address{Department of Mathematics,
     University of Joensuu,
     PL 111, 80101 Joensuu, Finland}
\email{villesamuli.normi@joensuu.fi}

\date{\today}


\begin{abstract}
We use Lagrangian formalism and jet spaces to derive a computational model to simulate multibody dynamics with holonomic constraints. Our approach avoids the traditional problems of drift-off and spurious oscillations. Hence even long simulations remain physically relevant. We illustrate our method with several numerical examples.
\end{abstract}

\maketitle

\section{Introduction}
The equations of motion of rigid body systems were established in the 19th century; for example Appell \cite{appell} gives a very nice and complete classical treatment of the subject. In the precomputer age \emph{solving} a differential equation meant algebraic manipulation of the system such that the solution could be explicitly obtained. Hence in \cite{appell} a lot of space is devoted to various specific situations and/or coordinate systems in which a complete solution could be obtained. Of course explicit solutions are usually impossible to obtain and hence the attention shifted gradually to qualitative, i.e. geometric, analysis of the properties of the solution set. This point of view lead to a modern differential geometric formulation of classical mechanics which can be found in \cite{arnold}. However, the techniques developed in \cite{appell} and \cite{arnold}, while interesting in proper context, are not necessarily suitable or useful for numerical treatment of multibody systems.

One of the earliest attempts at the numerical simulation of the multibody systems  was Baumgarte's stabilization method \cite{baumgarte}. Subsequently there has been quite much interest in this topic and in \cite{amirouche} \cite{gdjbayo}  \cite{schwerin} \cite{shabana}  many different formulations are described. However, in spite of the extensive literature on the subject the simulation of constrained multibody systems remains a challenging problem and it appears that no definitive solution to the problem has been found.

The major difficulty in the simulations is how to respect the constraints in long simulations and at the same time to avoid spurious oscillations which occur quite often in various stabilization methods. In this paper we propose a new method, based on our earlier work in \cite{julk31}, \cite{julk36} and \cite{arti14}, which addresses these issues.

Our computational model considers differential equations in \emph{jet spaces}. Hence the differential geometry is essential in our approach; nevertheless our formulation differs from the standard geometric model given in \cite{arnold}. We use Lagrangian formalism to derive the equations of motion; this is more natural in jet space context than Hamiltonian formalism. Constraints and possible \emph{invariants} (like conservation of energy) are taken into account by restricting the dynamics to an appropriate submanifold of a jet space. As a consequence quite much of the computing time is spent on projecting intermediate points to this submanifold. On the other hand the drift off is completely avoided. Since there are no spurious stiffness or oscillations we can use an explicit method in time integration. We have adapted a well-known Runge--Kutta method by Dormand and Prince with the classical step size control to our context \cite{hanowa}.

We consider only holonomic constraints in this paper. In case of point masses these kind of systems were already used as examples in our earlier papers. However, in multibody dynamics a good numerical representation of the orientation of the rigid body is somewhat involved, and hence the generalization of our approach to this case required some work. Nonholonomic problems, while certainly interesting, are beyond the scope of the present article and we refer to \cite{bloch}, \cite{bullew}, \cite{krupkova}, and \cite{rabrhe} for more details on this subject.

The structure of the article is as follows. In Section \ref{sec:basic} we briefly recall some geometric notions which are needed in our computational model.
In section  \ref{sec:one_body} we discuss the relevant variational principles leading to
a Lagrangian formulation of, and derive, the equations of motion for
one rigid body with external forces acting on it.
We believe it is best to explain the main ideas in this simple context because extension to many body case is then, due to the holonomicity of the constraints, just a matter of establishing an appropriate notation.

In Section  \ref{sec:many_body} we extend everything to an arbitrary number of rigid bodies and introduce invariants associated to multibody systems. We also recall for completeness what happens in the planar case: the whole model becomes significantly simpler because the orientations are trivial in this case. In Section  \ref{sec:comp} we describe the actual numerical subtasks which are needed in our implementation. The description is rather brief because we can use extensively algorithms which are explained in more detail in \cite{arti14}. Then in Section  \ref{sec:numer} we present our numerical results and finally in Section  \ref{sec:concl} give some conclusions and indicate some directions for future work.

\textbf{Acknowledgements} We are grateful for useful and motivating
discussions with Aki Mikkola and Asko Rouvinen from the Laboratory of
Mechatronics and Virtual Engineering, Lappeenranta University of
Technology, Finland.

\section{Basic tools}
\label{sec:basic}

For more information on standard differential geometry we refer to \cite{spivak} and on jets to \cite{saunders}. Further details about calculus of variations can be found in \cite{giahil}. Differential equations in the jet space context are discussed in \cite{pommaret}, \cite{seiler} and \cite{tarkhanov}.

All maps and manifolds are assumed to be smooth, i.e. infinitely differentiable. All analysis is local, hence various maps and manifolds need to be smooth or defined only in some appropriate subsets. To simplify the notation these subsets are not indicated.

\subsection{Manifolds and jets}
The $\ell$'th differential of a map $f\,:\,\mathbb{R}^m\mapsto \mathbb{R}^k$ is denoted by $d^\ell f$ and its value at $p$ by $d^\ell f_p$. The subscript is sometimes omitted for simplicity, if the point $p$ is clear from the context. Let $M$ be a manifold. The tangent bundle of $M$ is denoted by $TM$, and the tangent space at $p\in M$ by $T_pM$.  A distribution $\mathcal{D}$ is a map that associates to each point $p\in M$ a subspace $\mathcal{D}_p$ of $T_pM$.

Let $M$ be a submanifold of $\mathbb{R}^n$. The objects defined on $M$ can be taken to be defined on $\mathbb{R}^n$ without writing explicitly the inclusion map. The inner product in $\mathbb{R}^n$ is denoted by $\langle\cdot,\cdot\rangle$ and the same notation will be used also for inner products in $T_pM$ and $T_p\mathbb{R}^n$ as usual. We may regard $T_pM$ as a subspace of $T_p\mathbb{R}^n$. The orthogonal complement of $T_pM$, the normal space, is denoted by $N_pM$.

Let $\pi\,:\, \mathcal{E}\rightarrow\mathcal{B}$ be a bundle and let $\pi^q\,:\,J_q(\mathcal{E})\rightarrow\mathcal{B}$ be the bundle of $q$-jets of $\mathcal{E}$. In the sequel we will mainly consider the case where $\mathcal{E}$ is the trivial bundle $\mathcal{E}=\mathbb{R}\times\mathbb{R}^n$ with the projection $\pi\,:\, \mathbb{R}\times\mathbb{R}^n\rightarrow\mathbb{R}$. The coordinates of $J_q(\mathcal{E})$ (called jet coordinates) are denoted by $(t,y^1,\dots,y^n,y_1^1,\dots,y^n_q)$.

A \emph{section} of the bundle is a map $y\,:\,\mathcal{B}\rightarrow\mathcal{E}$ such that $\pi\circ y=\mathsf{id}$. In case of the trivial bundle the section is simply the graph of the map. In jet geometry it is usually more convenient to use sections than maps. For simplicity of notation we will use the same letter for maps and sections, the intended meaning being clear from the context.

\subsection{Some matrix groups and algebras}
The special orthogonal group is
\[
   \mathbb{SO}(n)=\big\{ A\in\mathbb{R}^{n\times n}\,\big|\, A^TA=I\ ,\ \det(A)=1\big\}.
\]
Let us introduce a convenient representation of elements of $\mathbb{SO}(3)$. We need this in order to define the orientation of the rigid body. This representation is discussed in detail in \cite{rabrhe} where it is derived using quaternions. First let $S^3\subset\mathbb{R}^4$ be the unit sphere and let $\theta\in S^3$. Then we set
\begin{equation}
\begin{aligned}    R=&\tilde H H^T=
\begin{pmatrix}
      -\theta^1&\theta^0&-\theta^3&\theta^2\\
      -\theta^2&\theta^3&\theta^0&-\theta^1\\
      -\theta^3&-\theta^2&\theta^1&\theta^0
        \end{pmatrix}
   \begin{pmatrix}
      -\theta^1&\theta^0&\theta^3&-\theta^2\\
      -\theta^2&-\theta^3&\theta^0&\theta^1\\
      -\theta^3&\theta^2&-\theta^1&\theta^0
        \end{pmatrix}^T\\
=&   2\begin{pmatrix}
         (\theta^0)^2+(\theta^1)^2-\tfrac{1}{2}&\theta^1\theta^2-\theta^0\theta^3&
      \theta^1\theta^3+\theta^0\theta^2\\
  \theta^1\theta^2+\theta^0\theta^3& (\theta^0)^2+(\theta^2)^2-\tfrac{1}{2}&
      \theta^2\theta^3-\theta^0\theta^1\\
   \theta^1\theta^3-\theta^0\theta^2&\theta^2\theta^3+\theta^0\theta^1&
       (\theta^0)^2+(\theta^3)^2-\tfrac{1}{2}
        \end{pmatrix}.
\end{aligned}
\label{rotaatio}
\end{equation}
It is now straightforward to check that $R\in\mathbb{SO}(3)$. The parameters $\theta^i$ are called \emph{ Euler} parameters.\footnote{a.k.a. \emph{Rodrigues} or \emph{Cayley-Klein} parameters.}  Let us note that
\begin{itemize}
\item[-] the rows of $\tilde H$ and $H$ are an orthonormal basis of $T_\theta S^3$, i.e. $HH^T=\tilde H \tilde{H}^T =I$.
\item[-] the parameters $\theta$ are not coordinates of $\mathbb{SO}(3)$ in the differential geometric sense.
\item[-] $\theta$ and $-\theta$ correspond to the same element $R$. Hence $S^3$ is a two sheeted covering space of $\mathbb{SO}(3)$ and consequently $\mathbb{SO}(3)$ is diffeomorphic to real projective space $\mathbb{RP}^3$.
\end{itemize}
Let us list some properties of $\tilde H$ and $H$ which are needed later. It is sometimes convenient to regard $\tilde H$ and $H$ as linear maps $\mathbb{R}^4\rightarrow\mathbb{R}^{3\times 4}$. In this way it is natural to write $H_1=H(\theta_1)$. The following formulas are easily verified by direct computation.
\begin{lem}
\begin{align*}
        & \tilde H^T\tilde H= H^T H=I-\theta \theta^T     \\
        & \tilde H_1H^T-\tilde H H_1^T =
          H_1\tilde H^T-H \tilde H_1^T=0 \\
 & H_1H^T+H H_1^T=  \tilde H_1\tilde H^T+\tilde H \tilde H_1^T =0.
\end{align*}
Moreover if $v$, $w\in\mathbb{R}^4$ are any vectors, then
\[
   H(v)w+H(w)v=\tilde H(v)w+\tilde H(w)v=0.
\]
In particular
\[
         \tilde H\theta=H\theta=\tilde H_1\theta_1=H_1\theta_1=0.
\]
\label{hmatriisi}
\end{lem}

From the differential geometric point of view $\mathbb{SO}(n)$ is a \emph{Lie group}, and the corresponding \emph{Lie algebra} is
\[
  \mathfrak{so}(n)=\big\{ A\in\mathbb{R}^{n\times n}\,\big|\, A^T=-A\big\}.
\]
Geometrically we may identify $\mathfrak{so}(n)$ with $T_I\mathbb{SO}(n)$. The important point for us is that $\mathfrak{so}(3)$ and $\mathbb{R}^3$ are naturally isomorphic as Lie algebras. To see this note that any $\hat\Omega\in\mathfrak{so}(3)$ is of the form
\[
  \hat\Omega =\begin{pmatrix}
           0&-\Omega^3&\Omega^2\\
            \Omega^3&0&-\Omega^1\\
            -\Omega^2&\Omega^1&0
            \end{pmatrix}.
\]
Hence putting $\Omega=(\Omega^1,\Omega^2,\Omega^3)$ we have for any $v\in\mathbb{R}^3$
\begin{equation}
        \hat\Omega v=\Omega\times v.
\label{liealg}
\end{equation}
One consequence of this isomorphism is that the tangent bundle of $\mathbb{SO}(3)$ is trivial:
\[
   T\mathbb{SO}(3)\simeq \mathbb{SO}(3)\times\mathbb{R}^3.
\]
From this it follows that the tangent bundle of $S^3$ is also trivial: $TS^3\simeq S^3\times\mathbb{R}^3$. This is very important from the computational point of view as we will see later.

\subsection{Variational problems}
Let $(M,G)$ be an $n$--dimensional Riemannian manifold. We denote the coordinates of $M$ by $y$. By the standard abuse of notation the same letter is also used for the curve $y\,:\, \mathbb{R}\rightarrow M$ and its coordinate representation. Then let us consider the variational problem where we want to find the extremals of the following map.
\begin{equation}
     J(y)=\int L(t,y, y_1) dt
\label{variint}
\end{equation}
where the integration interval, as well as the values of $y$ at the end points, are fixed during the variation.
The integrand is called the \emph{Lagrangian}.
\begin{defi} The Euler--Lagrange operator $\mathbb{E}_L$ and the Euler--Lagrange equations for the problem \eqref{variint} are
\begin{equation}
 \mathbb{E}_L(y)=     \frac{d}{dt}\frac{\partial L}{\partial y_1}-
        \frac{\partial L}{\partial y}=0.
\label{eulaeq}
\end{equation}
Extremals of \eqref{variint} are solutions of the Euler--Lagrange equations.
\end{defi}
Then given a Riemannian metric $G$ we may consider
\[
     J(y)=\tfrac{1}{2}\int \langle y_1, Gy_1\rangle dt.
\]
The extremals of this problem are called \emph{geodesics}. Let us recall the Euler--Lagrange equations for geodesics. First we set
\[
    [i\,j,k] = \frac{1}{2}\,\Big( \frac{\partial g_{ik}}{\partial y^j}+
    \frac{\partial g_{jk}}{\partial y^i}-\frac{\partial g_{ij}}{\partial y^k}\Big)
      \hspace{10mm}  i,j,k=1,\dots,n.
\]
These are called the \emph{Christoffel symbols} of the first kind. To write down the equations in a compact way we now introduce some nonstandard notation. Let us define the matrices
\[
    \mathsf{Chr}_k\in\mathbb{R}^{n\times n}\hspace{10mm}
     \big(\mathsf{Chr}_k\big)_{ij}=[i\,j,k].
\]
Note that the matrices $\mathsf{Chr}_k$ are symmetric. Putting all these matrices together we obtain a ``three dimensional'' object:
\[
    \mathsf{Chr}=\big(\mathsf{Chr}_1,\dots,\mathsf{Chr}_n\big)
       \in\mathbb{R}^{n\times n\times n}  \hspace{10mm}
       \mathsf{Chr}_{ijk}=[i\,j,k].
\]
Note that $\mathsf{Chr}$ is not a tensor: it transforms differently in coordinate changes. Anyway we can now view $\mathsf{Chr}$ as a bilinear map: given $v$, $w\in\mathbb{R}^n$ we can now set
\[
   \mathsf{Chr}(v,w)=\Big( \langle v,\mathsf{Chr}_1 w\rangle,\dots,
                         \langle v,\mathsf{Chr}_n w\rangle\Big).
\]
Note that $\mathsf{Chr}(v,w)=\mathsf{Chr}(w,v)$. Now a curve $y\,:\,\mathbb{R}\rightarrow M$ is a geodesic, if
\begin{equation}
   G y_2+\mathsf{Chr}(y_1,y_1)=0.
\label{geo}
\end{equation}
In differential geometry books the Christoffel symbols of the second kind are more often used than the symbols of the first kind. The symbols of the second kind are the components of the ``matrix'' $G^{-1}\mathsf{Chr}$. However, in computations one wants to avoid inverting matrices because this is both inefficient and unnecessary.

The equations for geodesics can be formulated in a coordinate free way using the (Levi-Civita) \emph{connection}. From this point of view the Christoffel symbols represent the connection in a coordinate system.

\subsection{Differential equations}
Let $\pi\,:\, \mathcal{E}\rightarrow\mathcal{B}$ be a bundle.
\begin{defi}
A (partial) differential equation of order $q$ on $\mathcal{E}$ is a submanifold $\mathcal{R}_q$ of $J_q(\mathcal{E})$.
\label{difyht}
\end{defi}

In jet coordinates the manifold $\mathcal{R}_q$ can be represented as a zero set of some map $f\,:\,J_q(\mathcal{E})\simeq\mathbb{R}^{(q+1)n+1}\rightarrow\mathbb{R}^k$:
\begin{equation}
\mathcal{R}_q\quad:\quad     f(t,y,y_1,\dots,y_q)=0.
\label{systeemi}
\end{equation}
To define solutions we introduce the following one forms
\begin{equation}
     \alfa_j^i=dy_{j-1}^i-y_j^idt
       \hspace{10mm}i=1,\dots,n
        \hspace{10mm}j=1,\dots,q.
\label{pfaff}
\end{equation}
Let $p\in J_q(\mathcal{E})$ and $v_p\in T_pJ_q(\mathcal{E})$ and let us further define distributions $\mathcal{C}$ and $\mathcal{D}$ by
\begin{equation}
\begin{aligned}
   \mathcal{C}_p=&\big\{v_p\in T_pJ_q(\mathcal{E})\,\big|\,
        \alfa_j^i(v_p)=0\big\},\\[2mm]
   \mathcal{D}_p=&T_p\mathcal{R}_q\cap C_p,
\end{aligned}
\label{distr}
\end{equation}
$\mathcal{C}$ is called \emph{Cartan distribution}. Now we can define the solutions as follows.
\begin{defi} Let $\mathcal{R}_q\subset J_q(\mathcal{E})$ be involutive and suppose that the distribution  $\mathcal{D}$ defined in \eqref{distr} is one-dimensional. A solution of $\mathcal{R}_q$ is an integral manifold of $\mathcal{D}$.
\label{ratk}
\end{defi}
The importance of involutivity from the point of view of numerical computations is discussed in \cite{julk31}. Intuitively a system is involutive if we cannot get new equations of order $q$ or less by differentiating the equations and eliminating the higher derivatives.

Now the equations of motion in Lagrangian mechanics take a very particular form and it turns out that using directly the above formulation of the problem would be very inefficient. Fortunately it is not difficult to adapt the jet space approach to this class of problems as we will see. Let us here briefly indicate what kind of changes we have to make to the general framework.

As is well known the Lagrangian dynamics with holonomic constraints gives a system of second order equations which can be written as
\begin{equation}
     f_\lambda(t,y,y_1,y_2,\lambda)=0
\label{e-l-system}
\end{equation}
where $y\,:\,\mathbb{R}\rightarrow\mathbb{R}^n$ gives the coordinates of the configuration space and $\lambda\,:\,\mathbb{R}\rightarrow\mathbb{R}^\ell$ is the Lagrange multiplier. Now if we regard $\lambda$ simply as any other dependent variable, then we would have to work with the space $J_2(\mathcal{E}_\lambda)$ where $\mathcal{E}_\lambda=\mathbb{R}\times\mathbb{R}^n\times\mathbb{R}^\ell$ and whose dimension is
\[
    \textrm{dim}\big(J_2(\mathcal{E}_\lambda)\big)=3(n+\ell)+1.
\]
Moreover the system \eqref{e-l-system} is not in the involutive form; in fact usually one has to differentiate and eliminate 4 times before we reach the involutive form, see \cite{julk31} for an example. However, it is unnecessary to treat $y$ and $\lambda$ in the same way, and we can work with the space $J_1(\mathcal{E})$ whose dimension is just
\[
    \textrm{dim}\big(J_1(\mathcal{E})\big)=2n+1.
\]
So our problem is geometrically as follows:
\begin{itemize}
\item the set of possible states of the system is given by a manifold $M\subset J_1(\mathcal{E})$
\item the dynamics of the system is given by a one dimensional distribution $\mathcal{D}$ on $M$
\end{itemize}
Before a detailed description of the computational model let us here outline the basic strategy. First the manifold $M$ is given as a zeroset of some map, and solutions are curves on $M$. To choose the right curve through some point $p$ we need to compute the appropriate distribution at $p$. Now given $p=(t,y,y_1)\in M$ it is easy to see that $v_p=(1,y_1,w)\in \mathcal{C}_p$ for any $w$. Hence our task is to compute the correct $w$. Next consider the section $j^1(y)\,:\,\mathbb{R}\rightarrow M\subset J_1(\mathcal{E})$ given by $t\mapsto (t,y(t),y'(t))$; the tangent vector to this curve is $(1,y'(t),y''(t))$, and hence it seems that we should take $w=y_2$. However, $y_2$ is ``outside'' of $J_1(\mathcal{E})$, so this is not directly applicable. But then we show that we can compute the correct $y_2$ with help of the Euler--Lagrange equations, and use this to define the right distribution.

\section{One rigid body}
\label{sec:one_body}
It turns out that all the relevant ideas that we will need in our computational model can be explained already in case of one rigid body, and extension to the general case is mainly a matter of introducing appropriate notation. Hence we feel that one gets a clearer picture of our approach by treating thoroughly the simple case and then rapidly indicating how to pass to the general case.

Equations of motion for rigid bodies can be found in most textbooks on classical mechanics. However, many formulations are not at all suitable for numerical computations. There are also many books which are devoted to computational aspects of multibody systems, see for example \cite{amirouche} \cite{gdjbayo} \cite{rabrhe} \cite{schwerin} \cite{shabana}. Since our model is new we cannot, however, always refer to standard literature for details, and hence we have to spend some time to describe our approach. In particular as far as we know jet spaces have not been previously used in multibody simulations.

\subsection{Configuration space and state space}
To describe the motion of a rigid body we need a fixed coordinate system (or spatial coordinate system), and the coordinate system moving with the body (or body coordinate system). The origin of the moving coordinate system is assumed to be at the centre of mass of the body. A typical point in spatial coordinates is denoted by $x$ and in body coordinates by $\chi$. The body itself is denoted by $\mathcal{B}\subset\mathbb{R}^3$, its mass density by $\rho$ and its mass by $m=m(\mathcal{B})$.
\begin{defi}
Let $v\in\mathbb{R}^3$ and set
\[
   \mathbb{I}v=\int_\mathcal{B} \big(\chi\times(v\times\chi)\big)\rho(\chi)d \chi
\]
$\mathbb{I}$ is the inertia tensor of body $\mathcal{B}$.
\end{defi}
One can readily check that $\mathbb{I}$ is symmetric and positive definite. The position of the centre of mass in the fixed coordinate system is given the map $r=(r^1,r^2,r^3)\,:\,\mathbb{R}\rightarrow\mathbb{R}^3$. To describe the motion of other points of the body we need to specify the orientation of the rigid body; this can be done with a rotation matrix.
Now the relation between spatial and body coordinates is given by the formula
\[
        x(t)=r(t)+R(t)\chi
\]
where $R(t)\in \mathbb{SO}(3)$. Hence the configuration space of one rigid body is $\mathbb{R}^3\times \mathbb{SO}(3)$. However, it is computationally not easy to work directly with the space $\mathbb{SO}(3)$. Instead we will choose the following framework.
\begin{defi} The configuration (resp. state) space of one rigid body is the bundle
\[
    \mathcal{Q}=\mathbb{R}\times\mathbb{R}^3\times S^3\rightarrow\mathbb{R}\hspace{10mm} \big(\ \ resp.\quad J_1(\mathcal{Q})\quad\big).
\]
\end{defi}
Geometrically this means that we replace $\mathbb{SO}(3)$ by its covering space $S^3$, using the representation given by \eqref{rotaatio}.  Recall that in jet context it is natural to think in terms of bundles; classically one thinks in terms of fibers, and hence the state space is the tangent bundle of the configuration space. This is really almost the same as our formulation because of the following identifications:
\[
  J_1(\mathcal{Q})\simeq \mathbb{R}\times T \big(\mathbb{R}^3\times S^3\big)
      \simeq  \mathbb{R}\times T \mathbb{R}^3\times TS^3.
\]
Let us then introduce the bundle
\[
\pi\,:\,\mathcal{E}=\mathbb{R}\times\mathbb{R}^7\rightarrow\mathbb{R}
\hspace{15mm}    y=(y^1,\dots,y^7)=(r,\theta).
\]
In this way we consider $\mathcal{Q}$ as a subset of $\mathcal{E}$. Then we define
\begin{align*}
   & f_\theta\,:\,J_1(\mathcal{E})\rightarrow\mathbb{R}^2\hspace{10mm}
       f_\theta(t,y,y_1)=\begin{pmatrix}
                 |\theta|^2-1\\
                 \langle \theta,\theta_1\rangle
                \end{pmatrix}\\
    &  M_\theta=f^{-1}(0)\subset J_1(\mathcal{E}).
\end{align*}

\begin{defi}
$M_\theta$ is the computational state space of one rigid body.
\end{defi}
This representation of the state space is much more convenient than the representation obtained by introducing coordinates on $S^3$ or $\mathbb{SO}(3)$. Neither of these spaces is diffeomorphic to $\mathbb{R}^3$, hence any coordinate system is necessarily local. So in general one should change coordinates during the computation, and this would be very annoying. On the other hand the parameters $\theta$, while not coordinates, provide anyway a global representation of the configuration space.

Finally let us note that $M_\theta$, being a submanifold of $J_1(\mathcal{E})$, is a differential equation according to Definition \ref{difyht}. However, it is obviously \emph{underdetermined}: the dimension of the distribution defined in \eqref{distr} is greater than one. Note that $\textrm{dim}\big(\mathcal{C}_p\big)=8$ and it is straightforward to check that
\[
   \textrm{dim}\big(T_p M_\theta\cap \mathcal{C}_p\big)=7.
\]
This dimension reflects in a natural way the degree of indeterminacy of the system: there are 6 degrees of freedom and the time variable.

So we have now an appropriate state space for the computational purposes. The next task is to introduce dynamics; i.e. choose an appropriate distribution on $M_\theta$.

\subsection{Variational principle}
The equations of motion of mechanical systems are usually derived from some variational principle. There are many essentially equivalent principles, see for example \cite[Tome 2, p. 451 and p. 492]{appell}, \cite[Chapter 3]{arnold} and \cite[p. 205]{bullew} for some discussion. In addition one has to choose between Hamiltonian and Lagrangian formalism. We use the latter because it is more natural in jet space context. The formulations of variational principles below are adapted from  \cite[Chapter 4]{bullew}. The statements are given geometrically, i.e. in a coordinate free way. Expressing them in a coordinate system give formulas which are found in classical textbooks.

We will use a variational principle to determine the distribution on $M_\theta$. To do this we first formulate the variational principle in $J_1(\mathcal{Q})$ and then interprete the result in terms of $M_\theta\subset J_1(\mathcal{E})$. The starting point is the kinetic energy of the rigid body. To define this we need to introduce \emph{angular velocities}.

Let $R\,:\,\mathbb{R}\rightarrow\mathbb{SO}(3)$ and define
\[
     \hat\Omega=R^TR_1 \hspace{15mm}
      \hat\omega=R_1R^T.
\]
These matrices belong to $\mathfrak{so}(3)$; hence we can define corresponding $\Omega$ and $\omega$ as in \eqref{liealg}. $\hat\Omega$ and $\Omega$ (resp. $\hat\omega$ and $\omega$) is called the \emph{ body} (resp. \emph{spatial}) \emph{angular velocity}.

\begin{defi}
The kinetic energy of a rigid body is given by
\begin{equation}
   T=T_{\mathsf{tr}}+T_{\mathsf{ro}}=
  \tfrac{1}{2}\,m|r_1|^2+ \tfrac{1}{2}\,\langle \Omega,\mathbb{I}\Omega\rangle
\label{kinen}
\end{equation}
where $T_{\mathsf{tr}}$ is the translational energy and $T_{\mathsf{ro}}$ the rotational energy.
\end{defi}
The angular velocities are related to $\theta_1$ by the following simple formulas.
\begin{lem}
\[
     \Omega=2H\theta_1 \hspace{15mm}  \omega=R\Omega=2\tilde H\theta_1.
\]
\end{lem}
\begin{proof}
Using Lemma \ref{hmatriisi} we compute
\[
   \hat \Omega=R^TR_1=H\tilde H^T(\tilde H_1 H^T+\tilde H H_1^T)=2H H_1^T
\]
and it is straightforward to check that $\hat{H \theta_1} =H H_1^T$, proving
the claim for $\Omega$.
Similar computation proves the formula for $\omega$.
\end{proof}
Hence the kinetic energy can be written as
\[
   T  =\tfrac{1}{2}\,m|r_1|^2+ 2\,\langle H\theta_1,\mathbb{I}H\theta_1\rangle.
\]

  Now the variational principle\footnote{a.k.a. the Jacobi principle \cite{lanczos}} we need is:
\begin{itemize}
\item[] \emph{The kinetic energy defines a Riemannian metric on the configuration space. The motions of the rigid body in the absence of forces are geodesics for this metric.}
\end{itemize}
The technical difficulty in using this principle is related to the rotational energy, so let us first treat the case that the rigid body is fixed at the centre of mass. The problem is that we cannot directly use the equations \eqref{eulaeq} or formulas \eqref{geo} because the parameters $\theta$ are not coordinates on $S^3$. Hence we have to introduce some local coordinates, and express the parameters $\theta$ in terms of these coordinates.

So let us choose a local parametrisation of $S^3$, i.e. we choose some open sets $U_1\subset\mathbb{R}^3$ and $U_2\subset S^3$ and a diffeomorphism $\varphi\,:\,U_1\rightarrow U_2$, $\theta=\varphi(\alpha)$. Hence we can write the rotational kinetic energy as
\begin{equation}
     T_{\mathsf{ro}}= 2\,\langle Hd\varphi\,\alpha_1,\mathbb{I}Hd\varphi\, \alpha_1\rangle
= \tfrac{1}{2}\,\langle \alpha_1,G_{\mathsf{ro}} \alpha_1\rangle
\label{roten}
\end{equation}
where the Riemannian metric $G_{\mathsf{ro}}$ is given by
\[
      G_{\mathsf{ro}}=4 d\varphi^TH^T \mathbb{I}Hd\varphi
\]
Now we could use the Euler--Lagrange equations \eqref{eulaeq} or formula \eqref{geo} to compute the equations of motion in terms of $\alpha$. However, we want to interprete the resulting equations in terms of $\theta$. It seems that in this case it is best to start again from Euler--Lagrange equations, rather than trying to express Christoffel symbols in terms of $\theta$.

Hence we need to compute the Euler--Lagrange operator for the Lagrangian $T_{\mathsf{ro}}$. First let us define
\[
   dG_{\mathsf{ro}}(v,w)=\Big( \big\langle v,\frac{\partial G_{\mathsf{ro}}}{\partial \alpha^1} w\big\rangle,
      \big\langle v,\frac{\partial G_{\mathsf{ro}}}{\partial \alpha^2} w\big\rangle,
       \big\langle v,\frac{\partial G_{\mathsf{ro}}}{\partial \alpha^3} w\big\rangle \Big)
\]
Thus $dG_{\mathsf{ro}}$ is a symmetric bilinear map $dG_{\mathsf{ro}}\,:\,\mathbb{R}^3\times \mathbb{R}^3\rightarrow\mathbb{R}^3$. Next we establish some formulas which are needed in the computations.
\begin{lem}
\begin{align*}
        d^2\varphi(\alpha_1,\cdot)&=\big(\frac{d}{dt} d\varphi\big)^T\\
        dH\,\theta_1&=-H_1\,d\varphi
\end{align*}
\label{kaavoja}
\end{lem}
\begin{proof} The first formula is quite straightforward to check. Then let us denote by $A_i$ (resp. $B_i$) the $i$th column of the left (resp. right) hand side of the second formula.Then using Lemma \ref{hmatriisi} we compute
\[
    A_i =\frac{\partial H}{\partial\alpha^i}\theta_1=
        H\big(\frac{\partial \varphi}{\partial\alpha^i}\big)\theta_1=
        -H(\theta_1)\frac{\partial \varphi}{\partial\alpha^i}=
           -H_1 \frac{\partial \varphi}{\partial\alpha^i}=-B_i
\]
\end{proof}

\begin{lem} The Euler--Lagrange operator for the rotational energy $T_{\mathsf{ro}}$ is
\begin{equation}
   \mathbb{E}_{T_{\mathsf{ro}}}(\alpha)=
 \frac{d}{dt}\frac{\partial T_{\mathsf{ro}}}{\partial \alpha_1}-
        \frac{\partial T_{\mathsf{ro}}}{\partial \alpha}=     4\,  d\varphi^TH^T \mathbb{I}H \theta_2+8\,  d\varphi^TH_1^T \mathbb{I}H\theta_1
\label{eularot}
\end{equation}
\label{eularotlem}
\end{lem}
\begin{proof}
First we compute
\begin{align*}
  &   \frac{d}{dt}\frac{\partial T_{\mathsf{ro}}}{\partial \alpha_1}=
     \frac{d}{dt}G_{\mathsf{ro}} \alpha_1=\\
   & 4\, d\varphi^TH^T \mathbb{I}H \frac{d}{dt}\big( d\varphi\, \alpha_1\big)+
     4\,  d\varphi^TH^T \mathbb{I}H_1\, d\varphi\,\alpha_1+\\
&   4\,  d\varphi^TH_1^T \mathbb{I}Hd\varphi\,\alpha_1+
     4\,  \big(\frac{d}{dt} d\varphi\big)^TH^T \mathbb{I}Hd\varphi\,\alpha_1=\\
   & 4\,  d\varphi^TH^T \mathbb{I}H \theta_2+4\,  d\varphi^TH_1^T \mathbb{I}H\theta_1+ 4\,  \big(\frac{d}{dt} d\varphi\big)^TH^T \mathbb{I}H \theta_1.
\end{align*}
Further we obtain
\[
   \big\langle \alpha_1,\frac{\partial G_{\mathsf{ro}}}{\partial \alpha^i} \alpha_1\big\rangle=
    8 \big\langle \frac{\partial d\varphi}{\partial \alpha^i} \alpha_1,H^T \mathbb{I} H\theta_1\big\rangle+
      8 \big\langle \frac{\partial H}{\partial \alpha^i} \theta_1,\mathbb{I} H\theta_1\big\rangle.
\]
Hence we can write
\[
   \frac{\partial T_{\mathsf{ro}}}{\partial \alpha}=
    \tfrac{1}{2} dG_{\mathsf{ro}}(\alpha_1,\alpha_1)=
    4 d^2\varphi(\alpha_1,\cdot) H^T \mathbb{I} H\theta_1+
      4 \big(dH\,\theta_1\big)^T\mathbb{I} H\theta_1.
\]
Then Lemma \ref{kaavoja} gives the result.
\end{proof}

\subsection{Forces and constraints}

In practice the rigid body is not free, but there are some forces acting on it. First let us introduce potential forces. Recall that the potential is simply a function $U$ in the configuration space. The variational principle can be extended to such systems in the following form:
\begin{itemize}
\item[] \emph{Motions of a rigid body fixed at the centre of mass subject to potential forces are given by the extremals of the Lagrangian }
\[
     L(\alpha,\alpha_1)=T_{\mathsf{ro}}-U(\alpha).
\]
\end{itemize}
Hence the computations in Lemma \ref{eularotlem} readily yield
\begin{koro} The Euler--Lagrange operator for $L=T_{\mathsf{ro}}-U$ is
\[
   \mathbb{E}_{L}(\alpha)=
  d\varphi^T\big(  4\, H^T \mathbb{I}H \theta_2+8\,  H_1^T \mathbb{I}H\theta_1+\nabla U\big).
\]
\end{koro}
Next we have the effect of external forces. Since the body is fixed at one point, we need only need to consider the torque acting on the body. Now given the torque $\tau$ in body coordinates, or $\tau_s=R\tau$ in spatial coordinates, how to incorporate this information in Euler--Lagrange equations? The appropriate concept in our context is given in the following definition. For a more general formulation we refer to \cite[p.188]{bullew}.
\begin{defi} $F_\tau$ is \emph{Lagrangian torque} corresponding to $\tau$, if
\[
   \langle F_\tau,\alpha_1\rangle=\langle \tau_s,\omega\rangle
\]
for all $\alpha_1$.
\end{defi}
The next version of the variational principle is as follows.
\begin{itemize}
\item[] \emph{Lagrangian torques are added to the right hand side of Euler--Lagrange equations.}
\end{itemize}
Then we have to compute $F_\tau$.
\begin{lem}
\[
    F_\tau=2d\varphi^T H^T\tau.
\]
\end{lem}
\begin{proof}
\begin{align*}
    \langle F_\tau,\alpha_1\rangle&=\langle \tau_s,\omega\rangle
                =\langle R\tau,2\tilde H\theta_1\rangle\\
        & =2 \langle\tilde H^T \tilde H H^T \tau,d\varphi\,\alpha_1\rangle
       =  2 \langle d\varphi^T H^T \tau,\alpha_1\rangle.
\end{align*}
\end{proof}
Hence our system so far can be written as
\[
   \mathbb{E}_{L}(\alpha)=
  d\varphi^T\big(  4\, H^T \mathbb{I}H \theta_2+8\,  H_1^T \mathbb{I}H\theta_1+\nabla U\big)= 2d\varphi^T H^T\tau.
\]
Let us then introduce holonomic constraints. We will suppose that the constraints are scleronomic, i.e. do not depend on time. Hence the constraints are given by the map
\[
      g_\alpha\,:\,U_1\rightarrow\mathbb{R}^\ell
\]
where $U_1$ is the domain of $\varphi$. Let us also set $g_\alpha=g\circ\varphi$. Now the system is required to stay in the subset of the configuration space defined by the zero set of $g_\alpha$:
\[
     M_\alpha=g_\alpha^{-1}(0)\subset U_1.
\]
To force the system to stay in the correct manifold the variational principle is modified as follows:
\begin{itemize}
\item[] \emph{In the presence of holonomic constraints a (fictious) constraint force $F_c$, which is normal to the constraint manifold, is added to the right hand side of Euler--Lagrange equations.}
\end{itemize}
Hence at present we have
\[
   \mathbb{E}_{L}(\alpha)=
  d\varphi^T\big(  4\, H^T \mathbb{I}H \theta_2+8\,  H_1^T \mathbb{I}H\theta_1+\nabla U\big)= 2d\varphi^T H^T\tau+F_c.
\]
\begin{lem} We have
\begin{equation}
       4\,  \mathbb{I}H\theta_2 +Hdg^T\lambda+8\,  H H_1^T \mathbb{I}H\theta_1
          +H\nabla U=2\tau
\label{eulalop}
\end{equation}
where $\lambda$ is the Lagrange multiplier.
\end{lem}
\begin{proof}
The rows of $dg_\alpha$ span $N_\alpha M_\alpha$. Hence
\[
   \mathbb{E}_{L}(\alpha)=
  d\varphi^T\big(  4\, H^T \mathbb{I}H \theta_2+8\,  H_1^T \mathbb{I}H\theta_1+\nabla U\big)= 2d\varphi^T H^T\tau-dg_\alpha^T\lambda
\]
for some $\lambda$. But then $dg_\alpha=dg\,d\varphi$ implies
\[
    d\varphi^T\big(    4\,  H^T \mathbb{I}H \theta_2
      +dg^T\lambda+ 8\,  H_1^T \mathbb{I}H\theta_1
    +\nabla U  -2 H^T\tau \big)=0.
\]
Since the columns of $d\varphi$ span $T_\theta S^3$ this is equivalent to
\[
     4\,  H^T \mathbb{I}H \theta_2
      +dg^T\lambda+ 8\,  H_1^T \mathbb{I}H\theta_1
    +\nabla U  -2 H^T\tau
    \in N_\theta S^3
\]
whence pre-multiplying by $H$ gives the desired formula.
\end{proof}
For convenience let us give the following
\begin{defi} The bilinear map
\[
    K(\theta_1,\theta_1)=  H H_1^T \mathbb{I}H\theta_1
\]
is called the \emph{Christoffel map}.
\end{defi}
Note that $K$ contains the information of the Christoffel symbols for the metric defined by $G_{\mathsf{ro}}$.

Recall that our goal is to compute $\theta_2$. The previous Lemma shows that we also have to compute $\lambda$ in the process. Both can be determined as follows.

\begin{alg}{\sf \hspace{20mm}\\[-5mm]
\begin{itemize}
\item given $(t,\theta,\theta_1)$ do
\begin{itemize}
\item solve $c$ and $\lambda$ from the following system
\begin{equation}
    \begin{cases}
     4\,  \mathbb{I}c +Hdg^T\lambda=2\tau
        -8\,  K  -H \nabla U  \\
       dg\,H^T c= |\theta_1|^2 dg\,\theta-d^2g(\theta_1,\theta_1)
     \end{cases}
\label{algdistyht0}
\end{equation}
\item set $\theta_2=H^Tc-|\theta_1|^2\theta$
\end{itemize}
\end{itemize}  }
\label{algdist0}
\end{alg}
\begin{proof}
Since $\theta_2\in\mathbb{R}^4\simeq T_\theta S^3\oplus N_\theta S^3$ we have the corresponding decomposition
\[
       \theta_2=H^Tc+a\theta
\]
for some $c$ and $a$. Differentiating the constraint $|\theta|^2-1=0$ twice it is seen that $a=-|\theta_1|^2$. Similarly differentiating the constraint $g$ twice we get
\[
    dg\,\theta_2+d^2 g(\theta_1,\theta_1)=0
\]
But then substituting the decomposition of $\theta_2$ to this equation (resp. to \eqref{eulalop})  gives the second (resp. first) equation in \eqref{algdistyht0}.
\end{proof}

Now that the hard work has been done the rest follows rather painlessly. If the centre of mass of the rigid body is not fixed, then the relevant Riemannian metric is
\[
     G=\begin{pmatrix}
        m\,I_3&0\\
         0& G_{\mathsf{ro}}
        \end{pmatrix}.
\]
Computing the appropriate Euler--Lagrange operator is straightforward. Then we have to also include the resultant of the external forces, denoted by $F$, which is applied at the centre of mass.  Further we set $d=(r_2,c)$, $F_e=(F,2\tau)$ (where $\tau$ is the torque in body coordinates as before) and
\begin{align*}
   \mathsf{H}&=\begin{pmatrix}
                  I_3&0\\ 0&H
                      \end{pmatrix} \in\mathbb{R}^{6\times 7}\hspace{10mm}
   \mathsf{I}=\begin{pmatrix}
             0&0\\  0 & |\theta_1|^2 I_4
              \end{pmatrix} \in\mathbb{R}^{7\times 7}\\[3mm]
   E&=\begin{pmatrix}
           m I_3&0\\
            0&4\,  \mathbb{I}
      \end{pmatrix}  \hspace{15mm}
       \tilde K=
           \begin{pmatrix}
              0\\
         K(\theta_1,\theta_1)
                 \end{pmatrix}.
\end{align*}
Next we introduce constraints. Let $g\,:\,\mathcal{E}\rightarrow\mathbb{R}^\ell$ be a map which does not depend on time and set
\[
        \mathbb{R}\times M_g=g^{-1}(0)\subset\mathcal{E}.
\]
Further we define
\begin{align*}
   &      f_{\mathsf{hc}}\,:\,  J_1(\mathcal{E})\rightarrow\mathbb{R}^{2\ell}
             \hspace{10mm}
           f_{\mathsf{hc}}(t,y,y_1)=\begin{pmatrix}
                       dg\,y_1\\ g
                    \end{pmatrix}\\
   &          M_{\mathsf{hc}}=f_{\mathsf{hc}}^{-1}(0)\subset J_1(\mathcal{E})
\end{align*}
where the subscript {\sf hc} refers to ``holonomic constraints''.
\begin{defi} The constrained configuration and state spaces are
\begin{align*}
    \mathcal{Q}_{\mathsf{co}}&=\mathcal{Q}\cap \big(\mathbb{R}\times M_g\big)
        \subset\mathcal{E}\\
      M&=M_\theta\cap M_{\mathsf{hc}}\subset J_1(\mathcal{E})
\end{align*}
\label{rajtila}
\end{defi}
Now with these notations we can compute the relevant Euler--Lagrange equations in very much the same way as above. Since this extension is routine we merely state the final result.
\begin{teo} Let us consider the motion of a rigid body with the following data:
\begin{itemize}
\item[-] Lagrangian is $L=T-U$ where $T$ is given by \eqref{kinen} and $U$ is the potential.
\item[-] The body is subject to force $F_e=(F,2\tau)$.
\item[-] Constrained state space is given as in Definition \ref{rajtila}.
\end{itemize}
 Then the following algorithm computes the distribution on $M$.
\begin{alg}{\sf \hspace{20mm}\\[-5mm]
\begin{itemize}
\item given $p=(t,y,y_1)\in M$ do
\begin{itemize}
\item solve $(d,\lambda)$ from the following system
\begin{equation}
     \begin{cases}
     E\,d +\mathsf{H}dg^T\lambda=F_e -8\, \tilde K-\mathsf{H}\nabla U\\
      dg\,\mathsf{H}^Td=dg\,\mathsf{I}\,y-d^2g(y_1,y_1)
       \end{cases}
\label{algdist1yht}
\end{equation}
\item set $y_2=\mathsf{H}^Td-\mathsf{I}y$.
\item set $\mathcal{D}_p=\mathsf{span}\{(1,y_1,y_2)\}$.
\end{itemize}
\end{itemize}  }
\label{algdist1}
\end{alg}
\end{teo}
We summarize our conclusions as follows:
\begin{itemize}
\item[]\emph{ Algorithm \ref{algdist1} determines a one dimensional distribution on $M$. The integral manifolds of this distribution give the motions of one rigid body.}
\end{itemize}
Traditionally one speaks about equations of motion for the rigid body. In this way we could say that the description of $M$ and $\mathcal{D}$ are our equations of motion. Note that in this model an integral manifold can always be represented by a curve $y\,:\,\mathbb{R}\rightarrow\mathcal{Q}_{\mathsf{co}}$. This is due to the first component of the distribution being always positive. Hence we will refer to the solutions as curves whenever convenient.

\subsection{Energy}
Finally we examine how the energy $W=T+U$ evolves.
\begin{lem} If $y$ is a solution, then
\[
    \frac{d\,W}{dt}=\langle \mathsf{H}y_1,F_e\rangle=
           \langle r_1,F\rangle+2\langle H\theta_1,\tau\rangle
            =\langle r_1,F\rangle+\langle \Omega,\tau\rangle
\]
\label{energialemma}
\end{lem}
\begin{proof} The energy can be written as
\[
   W=  \tfrac{1}{2}\,\langle \mathsf{H}y_1,E\,\mathsf{H}y_1\rangle+U(y)
\]
Hence
\begin{align*}
    \frac{d\,W}{dt}&= \langle \mathsf{H}y_1,E\,\mathsf{H}y_2\rangle+
        \langle \mathsf{H}y_1,E\,\mathsf{H}_1y_1\rangle+
         \langle \nabla U,y_1\rangle  \\
&= \langle \mathsf{H}y_1,E\,\mathsf{H}y_2\rangle+
             \langle \nabla U,y_1\rangle
     = \langle y_1,\mathsf{H}^T E\,d\rangle+
             \langle \nabla U,y_1\rangle
\end{align*}
because $\mathsf{H}_1y_1=(0,H_1\theta_1)=(0,0)$. Then using \eqref{algdist1yht} we obtain
\[
   \mathsf{H}^T E\,d =\mathsf{H}^TF_e-\mathsf{H}^T\mathsf{H}dg^T\lambda -8\, \mathsf{H}^T\tilde K-\mathsf{H}^T\mathsf{H}\nabla U
\]
Now the result follows from the following computations.
\begin{align*}
  &   \langle y_1, \mathsf{H}^T\mathsf{H}dg^T\lambda\rangle=
   \langle y_1, dg^T\lambda\rangle=
   \langle dg\,y_1, \lambda\rangle=0\\
   &     \langle y_1, \mathsf{H}^T\mathsf{H}\nabla U\rangle=
         \langle y_1, \nabla U\rangle\\
  &  \langle y_1, \mathsf{H}^T\tilde K \rangle=
      \langle \theta_1, H^THH_1\mathbb{I}H\theta_1\rangle=
          \langle H_1 \theta_1, \mathbb{I}H\theta_1\rangle=0
\end{align*}
\end{proof}

\subsection{Classical formulation}
Classically the equations of motion for the rigid body are often written as follows. Let $F$ be the resultant of the external forces and let $\tau$ be the torque acting on the body. Then the equations of motion can be written as
\begin{equation}
    \begin{cases}
         mr_2=F\\
        \mathbb{I}\,\Omega_1+\Omega\times \mathbb{I}\,\Omega=\tau
     \end{cases}
\label{eqmotionclassical}
\end{equation}
The first equation is simply the Newton's second law, and the second equation is sometimes called Euler's equation, and hence the whole system is called Newton--Euler equations. Appell \cite[p. 48]{appell} refers to this system combined with the energy balance equation  as \emph{sept \'equations universelles}.

This way of viewing things is not so convenient when we have a system of several interacting bodies with constraints. One reason is that in Newton--Euler equations one has to consider internal as well as external forces. On the other hand in the variational formulation only external forces appear.

The fact that the first equation is of second order while the second one is of first order already indicates that the variables $r$ and $\Omega$ are not ``on the same level''. Physically this is clear since $r$ describes position while $\Omega$ describes (angular) velocity. In any case
it would be hard to avoid the conclusion that the Newton--Euler system is not very suitable for numerical computations, see \cite{rabrhe} for a detailed discussion.

\section{Systems of Rigid bodies}
\label{sec:many_body}

\subsection{Extension to general case}

Let us now formulate the problem with arbitrary number of rigid bodies. Let the number of bodies be $n_b$ and let the configuration space be
\[
    \mathcal{Q}=\mathbb{R}\times\Big( \mathbb{R}^3\times S^3\Big)^{n_b}
    \subset \mathcal{E}=\mathbb{R}\times\mathbb{R}^{7n_b}
\]
The coordinates of $\mathcal{E}$ are
\[
   (t,y)=\big(t,r^{(1)},\theta^{(1)},\dots,r^{(n_b)},\theta^{(n_b)}\big)
\]
where $(r^{(i)},\theta^{(i)})$ denote the position and Euler parameters of the $i$th body. The computational state space  $M_\theta\subset J_1(\mathcal{E})$ is the zero set of  the map
\begin{equation}
  f_\theta\,:\,J_1(\mathcal{E})\rightarrow\mathbb{R}^{2n_b}\hspace{10mm}
 f_\theta(t,y,y_1)=       \begin{pmatrix}
                  |\theta^{(1)}|^2-1\\[1mm]
                  \big\langle \theta^{(1)}, \theta_1^{(1)}\big\rangle\\
                     \vdots\\
                  |\theta^{(n_b)}|^2-1\\[1mm]
                  \big\langle \theta^{(n_b)}, \theta_1^{(n_b)}\big\rangle
                 \end{pmatrix}
 \label{perusmonisto}
\end{equation}

Let us further define
\begin{equation}
\begin{aligned}
      E&=\textrm{diag}\big(E^{(1)},\dots,E^{(n_b)} \big)
  & \mathsf{I}=\textrm{diag}\big(
                  \mathsf{I}^{(1)},\dots,\mathsf{I}^{(n_b)}\big)\\
   \mathsf{H}&=\textrm{diag}\big(
                  \mathsf{H}^{(1)},\dots,\mathsf{H}^{(n_b)}\big)
   &   F_e=\big(F_e^{(1)},\dots,F_e^{(n_b)}\big)\\
  d&=
\big(r_2^{(1)},c^{(1)},\dots,r_2^{(n_b)},c^{(n_b)}\big)
  &  \tilde K=\big(\tilde K^{(1)},\dots,\tilde K^{(n_b)}\big)
\end{aligned}
\label{yltapaus}
\end{equation}
where $E^{(i)}$ etc denotes the corresponding matrix or vector for the $i$th rigid body. The Lagrangian and the energy of the system  can now be written as
\begin{align*}
   L=T-U=&
   \sum_{i=1}^{n_b}
  \tfrac{1}{2}\,m_i|r_1^{(i)}|^2+ 2\,\langle H^{(i)}\theta_1^{(i)},\mathbb{I}^{(i)}H^{(i)}\theta_1^{(i)}\rangle-U(y)\\
   =&\tfrac{1}{2}\,\langle \mathsf{H}y_1,E\mathsf{H}y_1\rangle-U(y)\\
   W=T+U=&  \tfrac{1}{2}\,\langle \mathsf{H}y_1,E\,\mathsf{H}y_1\rangle+U(y).
\end{align*}
Next we introduce constraints in the same way as in the case of one rigid body, see Definition \ref{rajtila}. We set
\begin{equation}
\begin{aligned}
   & g\,:\,\mathcal{E}\rightarrow\mathbb{R}^\ell    &     
   \mathbb{R}\times M_g=g^{-1}(0)\subset\mathcal{E}  \\
   &      f_{\mathsf{hc}}\,:\, J_1(\mathcal{E})\rightarrow\mathbb{R}^{2\ell}
             &
           f_{\mathsf{hc}}(t,y,y_1)=\begin{pmatrix}
                       dg\,y_1\\ g
                    \end{pmatrix}\\
   &          M_{\mathsf{hc}}=f_{\mathsf{hc}}^{-1}(0)\subset J_1(\mathcal{E})\\
   & \mathcal{Q}_{\mathsf{co}}=\mathcal{Q}\cap \big(\mathbb{R}\times M_g\big)  &
      M=M_\theta\cap M_{\mathsf{hc}}.
\end{aligned}
\label{yltapaus2}
\end{equation}
Since we have $\ell$ independent constraints one says that there are $6n_b-\ell$ degrees of freedom in the system. In our context this can be checked by computing that
\[
   \textrm{dim}\big(T_pM\cap \mathcal{C}_p\big)=6n_b-\ell+1.
\]
Recall that ``$+1$'' is because of the time variable.

Then we get the following result. The proof is omitted because it is essentially the same as in the case of one rigid body.

\begin{teo} Adopting the notations in \eqref{yltapaus} and \eqref{yltapaus2}, Algorithm \ref{algdist1} computes the distribution also in the case of several rigid bodies.  Moreover we have
\[
    \frac{d\,W}{dt}=\langle \mathsf{H}y_1,F_e\rangle.
\]
\end{teo}
Hence in absence of external forces the energy remains constant. Moreover the constraint forces have no effect on energy. This is sometimes expressed by saying that constraint forces do no work.

Note that the interconnection forces between different bodies do not appear in this formulation; we need only specify the external forces acting on the system.
On the other hand if Newton--Euler equations are used, it is necessary to take care of interconnection forces as well. Hence it is by no means obvious that both models really yield the same results. However, a detailed discussion of the equivalence of variational and Newton--Euler models is beyond the scope of our article and we refer to \cite[Chapter 4]{bullew} for more information and references on this topic.

Now if we have external forces acting on the system it is convenient to define $W_{\mathsf{ext}}$, the work done by external forces:
\[
    W_{\mathsf{ext}}(t)= \int_0^t \langle \mathsf{H}y_1,F_e\rangle\,ds
\]
In this way the total energy
\[
   W_{\mathsf{total}}=T+U-W_{\mathsf{ext}}
\]
remains constant. This is useful in computations because we can compute $W_{\mathsf{ext}}$ approximatively using (for example) trapezoidal rule, and use this to control the total energy balance.

\subsection{Invariants}

In addition to constraints, there may be invariants\footnote{Invariants are also called constants of motion or sometimes even dynamical constraints. In the latter case our constraints are then kinematic constraints.} associated to the system which are modelled by a map $f_{\mathsf{inv}}\,:\,J_1(\mathcal{E})\rightarrow\mathbb{R}^r$. While invariants and constraints might superficially seem similar concepts, in reality they are of very different nature. Constraints are externally imposed on the system. Mathematically this means that Lagrange multipliers are needed in the corresponding equations. Physically this means that we introduce fictious forces which realize constraints. Invariants on the other hand are logical consequences of the dynamics. Hence their description requires neither Lagrange multipliers nor fictious forces.

A typical example of an invariant is the conservation of energy: in the absence of external forces we have seen that $W$ is an invariant. Usually the number of invariants is quite small, and indeed in many cases the energy is the only invariant. The invariants define a manifold
\[
     M_{\mathsf{inv}}=f_{\mathsf{inv}}^{-1}(0)\subset J_1(\mathcal{E}).
\]
Geometrically $M_{\mathsf{inv}}\cap M$ can be interpreted as a submanifold of $M$ which is invariant by the flow induced by the dynamics of the system. Hence the presence of invariants has no effect on the computation of the distribution.

\subsection{Planar case}

Of course the planar case is included in the above considerations as a special case. However, treating planar case directly is much more efficient computationally. Hence we briefly indicate the appropriate model.

Now in this case the relation between spatial and body coordinates is given by
\[
      x(t)=r(t)+R(t)\chi
\]
where $R\in \mathbb{SO}(2)$. But using the correspondence
\[
        \beta\mapsto \begin{pmatrix}
               \cos(\beta)&-\sin(\beta)\\
               \sin(\beta)&\cos(\beta)
                  \end{pmatrix}
\]
it is seen that the configuration space is $\mathbb{R}^2\times S^1$. Computationally we can regard $\beta\in\mathbb{R}$, with the understanding that the relevant data should be $2\pi$ periodic.  Next, it is evident that
\[
     \omega=\Omega=(0,0,\beta_1)
\]
and the inertia is described by the scalar
\[
      \mathbb{I}_p=\int_\mathcal{B}  |\chi|^2\rho(\chi)d\chi.
\]
The torque is also a scalar: $\tau\simeq (0,0,\tau)$. Then setting $y=(r,\beta)$, $F_e=(F,\tau)$, denoting the potential by $U$ and defining
\[
    E= \begin{pmatrix}
               m I_2&0\\
              0 &  \mathbb{I}_p
                  \end{pmatrix}
\]
we can write the Euler--Lagrange equations as
\[
    Ey_2+\nabla U=F_e.
\]
This formula can of course be generalised to the case of system of $n_b$ planar rigid bodies by following conventions:
\begin{align*}
      E&=\textrm{diag}\big(E^{(1)},\dots,E^{(n_b)} \big)\\
      F_e&=\big(F^{(1)},\tau^1,\dots,F^{(n_b)},\tau^{n_b}\big).
\end{align*}
Finally if the constraints are given by $g$ we see that the relevant distribution can be computed from the following equations.
\begin{equation}
   \begin{cases}
          Ey_2+dg^T\lambda+\nabla U=F_e\\
          dg\,y_2=-d^2g(y_1,y_1).
    \end{cases}
\label{planar}
\end{equation}
In particular it is seen that no term corresponding to Christoffel map appears in the planar case. The determination of the (constrained) state space is handled in the same way as in the general case.

\section{Computations}
\label{sec:comp}
The computational problem has 2 main ingredients:
\begin{itemize}
\item[(i)] given $p\in M$, compute $\mathcal{D}_p$
\item[(ii)] given $a\in J_1(\mathcal{E})$, project $a$ to $M$.
\end{itemize}
Of course in the latter problem $a$ must be sufficiently close to $M$ so that the projection is well defined. However, this is not a problem in practice since this condition is always satisfied if the step size is sufficiently small.

\subsection{Distribution}

How to solve the system \eqref{algdist1yht} as efficiently as possible? The following approach was already used in \cite{julk36} where more details and references can be found. Consider the block matrix
\[
      C  =\begin{pmatrix}
           A&B^T\\[2mm]
                  B&D
            \end{pmatrix}.
\]
If $A$ is invertible the \emph{Schur complement} of $A$ is
\[
         \mathsf{S}=D-BA^{-1}B^T.
\]
Now the Schur complement is useful in solving the linear system $Cw=b$ if
\begin{itemize}
\item[--] $A$ is ``easily invertible'' and
\item[--] $D$ (and hence $\mathsf{S}$) is (much) smaller than $A$.
\end{itemize}
But this is precisely the situation in our case! The system \eqref{algdist1yht} can be written as
\[
    \begin{pmatrix}
            E & \mathsf{H}(dg)^T\\[2mm]
                   dg\, \mathsf{H}^T &0
            \end{pmatrix}
\begin{pmatrix}
           d\\[2mm] \lambda
              \end{pmatrix}=
     \begin{pmatrix}
     F_e -  8 \tilde K-\mathsf{H}\nabla U \\[2mm]
           dg\,\mathsf{I}\,y-d^2g(y_1,y_1)
     \end{pmatrix}.
\]
The Schur complement is now
\[
         \mathsf{S}=- dg\, \mathsf{H}^T E^{-1}\mathsf{H}(dg)^T.
\]
Typically the size of $\mathsf{S}$ is much smaller than the size of $E$. But even more importantly, $E$ is indeed easily invertible, because it is block diagonal matrix with $6\times 6$ blocks. Hence we solve the system \eqref{algdist1yht} as follows:
\begin{alg}{\sf \hspace{20mm}\\[-5mm]
\begin{itemize}
\item solve $E u^1= F_e -  8 \tilde K-\mathsf{H}\nabla U$
\item compute $v^1= dg\, \mathsf{H}^T u^1-dg\,\mathsf{I}\,y+d^2g(y_1,y_1)$
\item solve $\mathsf{S}v^2=v^1$.
\item compute $u^2=\mathsf{H}(dg)^T v^2$
\item solve $ E u^3=u^2$
\item set $(d,\lambda)=\big(u^1+u^3,-v^2\big)$.
\end{itemize} }
\label{linsys}
\end{alg}
Note that the system $\mathsf{S}v^2=v^1$ must be solved iteratively while $E$ systems are solved ``exactly'' by a direct method. Then we get the following overall algorithm for the computation of the distribution.

\begin{alg}{\sf \hspace{20mm}\\[-5mm]
\begin{itemize}
\item given $p=(x,y,y_1)\in M$ do
\begin{itemize}
\item solve $(d,\lambda)$ from the system \eqref{algdist1yht} using Algorithm \ref{linsys}
\item set $y_2=\mathsf{H}^Td-\mathsf{I} y$
\item set $\mathcal{D}_p=\textrm{span}\{(1,y_1,y_2)\}$.
\end{itemize}
\end{itemize}  }
\label{algdistmanyconstr}
\end{alg}

\subsection{Projection}

We have in fact 3 relevant manifolds which we have to consider: $M_\theta$, $M_{\mathsf{hc}}$ and $M_{\mathsf{inv}}$. Let us first consider $M_\theta$. Since this involves only Euler parameters we formulate the projection directly in terms of $\theta$. Moreover different rigid bodies do not interact in the projection, so without loss of generality we consider just one rigid body.

In principle the orthogonal projection of $(a,a_1)\in T_{a}\mathbb{R}^4$ to $(\theta,\theta_1)\in T_\theta S^3$ could be computed by solving the system
\begin{equation}
 \begin{cases}
    \theta+\mu_1\theta+\mu_2\theta_1=a\\
     \theta_1+\mu_2\theta=a_1\\
    |\theta|^2-1=0\\
    \langle \theta,\theta_1\rangle=0.
  \end{cases}
\label{thetaproj}
\end{equation}
However, this would require Newton's method. Instead we project with the following algorithm.
\begin{alg}{\sf \hspace{20mm}\\[-5mm]
\begin{itemize}
\item given $(a,a_1)\in T_a\mathbb{R}^4$ do
\begin{itemize}
\item set $\theta=a/|a|$
\item set $\theta_1=a_1-\langle a_1,\theta\rangle \theta$.
\end{itemize}
\end{itemize}  }
\label{algproj1}
\end{alg}
Note that the point thus obtained satisfies 3 last equations of \eqref{thetaproj}. However, the first equation will be satisfied for some $\mu_1$ only if $\langle a,a_1\rangle=0$. Anyway in our application $\langle a,a_1\rangle$ will always be ``small'', more precisely $\langle a,a_1\rangle=O(h)$, so we may say that the projection is ``almost'' orthogonal.

Next let us consider $M_{\mathsf{hc}}$. Recall that we suppose that the constraints do not depend on $t$. Then given a point $(t,a,a_1) \in J_1(\mathcal{E})$, its orthogonal projection to $(t,y,y_1)\in M_{\mathsf{hc}}$ can be computed by solving the following system.
\begin{equation}
  F(y,y_1,\alpha,\beta)=\begin{pmatrix}
         y+  d^2g^T(y_1,\cdot)\alpha+dg^T\beta-a\\
        y_1+dg^T\alpha-a_1\\
        dg\, y_1\\g(y)
        \end{pmatrix}.
\label{newtonsysholo}
\end{equation}
But now we can do an ``almost orthogonal'' projection mimicking Algorithm \ref{algproj1} because the form of the equations is essentially the same in both cases. Given $(t,a,a_1)\in J_1(\mathcal{E})$ we can orthogonally project $a$ to $y\in M_g$ by solving the system
\begin{equation}
\begin{cases}
   y+(dg)^T\mu-a=0\\
      g(y)=0.
\end{cases}
\label{projconf}
\end{equation}
But having computed this we can solve $y_1$ and $\alpha$ from the system
\begin{equation}
\begin{cases}
   y_1+(dg)^T\alpha-a_1=0\\
   dg\,y_1=0.
\end{cases}
\label{projlin}
\end{equation}
Note that this is a linear system. Then we have obtained values $(y,y_1,\alpha)$ such that the last three equations of the system \eqref{newtonsysholo} are satisfied.

\begin{alg}{\sf \hspace{20mm}\\[-5mm]
\begin{itemize}
\item given $(t,a,a_1)\in J_1(\mathcal{E})$ do
\begin{itemize}
\item project $a$ to $y\in M_g$ by solving \eqref{projconf}
\item solve $y_1$ and $\alpha$ from system \eqref{projlin}.
\end{itemize}
\end{itemize}  }
\label{algprojholo}
\end{alg}

Finally we have the invariants, given by the map $f_{\mathsf{inv}}$. Now in general this map has no special structure, so given $\mathsf{a}=(t,a,a_1)\in J_1(\mathcal{E})$ we simply solve
\begin{equation}
\begin{cases}
   p+(df_{\mathsf{inv}})^T\mu-\mathsf{a}=0\\
      f_{\mathsf{inv}}(p)=0.
\end{cases}
\label{projinv}
\end{equation}

All these projections are then combined as follows.

\begin{alg}{\sf \hspace{20mm}\\[-5mm]
\begin{itemize}
\item given $p_{init}\in J_1(\mathcal{E})$ do
\begin{itemize}
\item set $p^0=p_{init}$, $j=0$
\item repeat
\begin{itemize}
\item[] project $p^j$ to $p_\theta\in M_\theta$ using Algorithm \ref{algproj1}
\item[] project $p_\theta$ to $p_{hc}\in M_{\mathsf{hc}}$ using Algorithm \ref{algprojholo}
\item[] project $p_{hc}$ to $p^{j+1}\in M_{\mathsf{inv}}$ by solving the system \eqref{projinv}
\end{itemize}
\item[] until convergence.
\end{itemize}
\end{itemize}  }
\label{algprojconstr}
\end{alg}

\subsection{Constraints in practice}

In practice the constraints for rigid bodies are of quite
particular form. In fact it turns out that
there are only three basic constraints and all needed constraints can be constructed by taking appropriate combinations of the basic ones.

The first one is a \emph{coincidence constraint} (or briefly \emph{C--constraint}), given by a map $g^{\mathsf{c}}:\mathbb{R}^{14}\to\mathbb{R}^3$. This constraint simply says that given points $\chi^{(i)}$ and $\chi^{(j)}$  in the coordinate systems of bodies $\mathcal{B}_i$ and $\mathcal{B}_j$ are really the same point in the spatial coordinate system. Hence
\begin{equation}
\label{g_C}
 g^{\mathsf{c}}(\mathcal{B}_i, \mathcal{B}_j) = r^{(j)} + R^{(j)}\chi^{(j)} - r^{(i)} - R^{(i)}\chi^{(i)} = 0.
\end{equation}
Mathematically $\chi^{(i)}$ and $\chi^{(j)}$ can be arbitrary, but in practice they are the positions of the relevant joint in the corresponding body coordinate systems.  We remind the reader that $r^{(i)}$ is the position of the centre of mass of the body $\mathcal{B}_i$ and  $R^{(i)}$ defines its orientation in the spatial (global) coordinate system.

Next we introduce a basic constraint where we require that two unit vectors $a^{(i)}, a^{(j)}$ given in body coordinate systems must be perpendicular to each other in the spatial coordinate system.
This is called \emph{symmetric orthogonality constraint} (or \emph{SO--constraint}),
and is given by the following map $g^{\mathsf{so}}:\mathbb{R}^8\to\mathbb{R}$:
\begin{equation}
\label{g_SO}
g^{\mathsf{so}}(a^{(i)}, a^{(j)}) = \langle R^{(i)} a^{(i)}, R^{(j)} a^{(j)}\rangle = 0.
\end{equation}

In the third constraint we are given a unit vector $ a^{(i)}$ and the points $\chi^{(i)}$ and $\chi^{(j)}$ in body coordinates. Let us consider the difference of $\chi^{(i)}$ and $\chi^{(j)}$ in spatial coordinates:
\[
   d^{(i,j)} = r^{(j)} + R^{(j)}\chi^{(j)} - r^{(i)} - R^{(i)}\chi^{(i)}.
\]
Now we require that this is orthogonal to $ a^{(i)}$ which must naturally also be expressed in the spatial coordinates. This is called \emph{nonsymmetric orthogonality constraint} ( or \emph{O--constraint}), and thus is given by  $g^{\mathsf{o}}:\mathbb{R}^{14}\to\mathbb{R}$:
\begin{equation}
\label{g_O}
\begin{aligned}
 g^{\mathsf{o}}(a^{(i)}, d^{(i,j)})& =  \langle R^{(i)} a^{(i)}, d^{(i,j)}\rangle \\
    & = \langle R^{(i)}a^{(i)}, r^{(j)} + R^{(j)}\chi^{(j)} - r^{(i)} \rangle - \langle a^{(i)}, \chi^{(i)}\rangle = 0.
\end{aligned}
\end{equation}
Note that an O--constraint has a singularity if $d^{(i,j)}$ happens to be zero. Table \ref{joints} indicates how one can define some typical joints using different combinations of the basic constraints.

\begin{table}[!htb]
\caption{Different types of joints.}\label{joints}
\begin{center}
\begin{tabular}{|l||c|c|}
\hline
type of joint & types of constraints & $\#$ of conditions\\
\hline
spherical      &
          1 C & 3 \\
\hline
 universal      &
  1 C and 1 SO & 4 \\
\hline
revolute       &
  1 C and 2 SO & 5 \\
\hline
cylindrical   &
  2 SO and 2 O & 4 \\
\hline
translational  &
  3 SO and 2 O  & 5\\
\hline
\end{tabular}
\end{center}
\end{table}

\section{numerical results}
\label{sec:numer}
\subsection{Implementation}
We use a 5th order Runge--Kutta scheme by Dormand and Prince \cite{hanowa}. In \cite{julk36} and \cite{arti14} we have explained in detail how this scheme (and other Runge--Kutta methods) can be adapted to jet space context. To speed up the computation
system \eqref{projinv} was solved for 4th and 5th order approximations only, not for the intermediate points. The tests indicated that this omission did not affect the quality of the computed solutions. To solve \eqref{projconf} and  \eqref{projinv} we used the \emph{inexact Newton} method \cite{deeist}; how to apply this in our context is explained in \cite{arti14}.

\subsection{Description of the models.}

In the following examples rigid bodies are homogeneous solids with mass density $\rho = 7810$.\footnote{We use standard SI units in our models.}
Gravitational acceleration is $\mathsf{g}=9.81$ in the negative direction of $x_2$-axis in spatial coordinates. The basic unit vectors are denoted by
\[
  \mathsf{e}_1 := (1,0,0) \qquad \mathsf{e}_2 := (0,1,0) \qquad \mathsf{e}_3 := (0,0,1).
\]
In Tables below we will specify the constraints by giving some vectors $a^{(i)}, a^{(j)}, b^{(i)}$ and $ b^{(j)}$ whose interpretation is as follows:
\begin{align*}
a^{(i)}  &: \hspace{1cm}\textrm{orthogonal to axis of a joint in $i$th body coordinates} \\
b^{(i)}  &: \hspace{1cm}\textrm{orthogonal to $a^{(i)}$ and $a^{(j)}$ in $i$th body coordinates} \\
b^{(j)}  &: \hspace{1cm}\textrm{parallel to $b^{(i)}$ in $j$th body coordinates}.
\end{align*}
Given these vectors we will need the following orthogonality constraints:
\begin{align*}
\textrm{SO}&: \qquad g^{\mathsf{so}}(a^{(i)}, a^{(j)}), ~g^{\mathsf{so}}(b^{(i)}, a^{(j)}), ~g^{\mathsf{so}}(a^{(i)}, b^{(j)}) \\
\textrm{O}&: \qquad  g^{\mathsf{o}}(a^{(i)}, d^{(i,j)}),~g^{\mathsf{o}}(b^{(i)}, d^{(i,j)}).
\end{align*}
If now in Table \ref{joints} we need a certain number of SO-- or O--constraints, then we take the corresponding constraints from the above list starting from left. For example the universal joint needs just $a^{(i)}$ and $a^{(j)}$ while only the translational joint uses $b^{(j)}$.

\subsubsection{3D pendulum}

This is a problem of three degrees of freedom. Our pendulum is
initially at rest along the spatial $x_1$-axis.  It is attached to a
spherical joint, which lies in the spatial origin.  There are no external forces acting on the system, except gravity; hence the energy $W=T+U$ remains constant.  We also tested with an imposed constant torque which showed an interesting behaviour and comment briefly on that.

\begin{figure}[h!]
\begin{center}
\centering
\epsfig{file=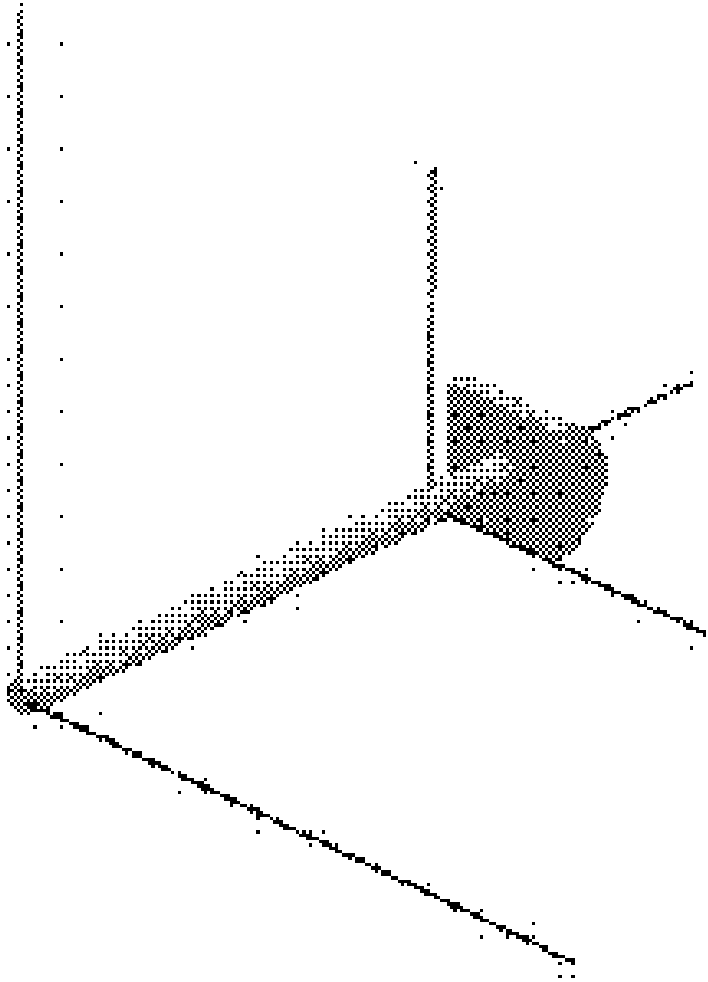,width=0.6\textwidth}
\caption{Initial configuration of the pendulum.}
\label{fig:heiluri_alkutila}
\end{center}
\end{figure}

Initial configuration is shown in Figure \ref{fig:heiluri_alkutila},
where the spatial coordinate axes are as follows: $\mathsf{e}_1$
towards the upper right corner, $\mathsf{e}_2$ upwards, $\mathsf{e}_3$
towards the lower right corner. We used the program \textsf{SolidWorks}\footnote{
\textsf{http://www.solidworks.com}.} to calculate the inertial tensor
\[
   \mathbb{I}=\textrm{diag}\big(0.147, 3.175, 3.154\big)
\]
as well as the centre of mass $m = 38.34$ of the pendulum.  Initially
the centre of mass in spatial coordinates is $r = \big(
0.765,0,0\big)$, and the orientation is given by Euler parameters
$\theta = (1,0,0,0)$. The vectors describing the position of the joint
are $\chi^{(0)} = (0,0,0)$ and $\chi^{(1)} = (-0.765,0,0)$.

\subsubsection{Planar quadrangle with joints}

\begin{figure}[h!]
\begin{center}
\centering
\epsfig{file=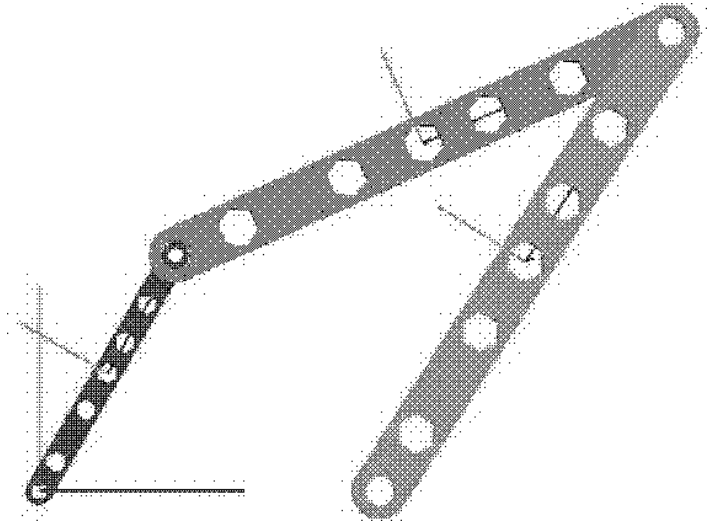,width=0.6\textwidth}
\caption{Initial configuration of the planar quadrangle.}
\label{fig:tasonivel_alkutila}
\end{center}
\end{figure}

This is a problem of one degree of freedom. Since this is a planar problem it is of course inefficient to use a 3D code to solve it. However, to model it we use 4 different types of joints, so this problem is very suitable for testing basic constraints.
The kinematic chain consists of three rigid bodies and a fixed ground body.
Each body, including the ground one, is attached to an adjacent one through a joint.
For testing purposes we have chosen revolute joints only at the ground body,
whereas the other two joints are spherical or cylindrical.
There is a single external torque acting on one of the bodies.
The spatial coordinate axes are in the Figure \ref{fig:tasonivel_alkutila}
as follows: $\mathsf{e}_1$
towards right,  $\mathsf{e}_2$ upwards,  $\mathsf{e}_3$
towards the reader.

In Table \ref{TNNK:table3}, where the corresponding setup is presented, symbol $\mathcal{B}_i$ (resp. $\mathcal{B}_j$)
stands for the ``first body'' (resp. ``second body'').
From the same table we find that the joints restrict 17 degrees of freedom from the system.
From Tables \ref{TNNK:table1} and \ref{TNNK:table2} one can find the chosen parameters of the model.

\begin{table}[!htb]
\caption{Parameters of the planar quadrangle.}\label{TNNK:table1}
\begin{center}
\begin{small}
\begin{tabular}{|c||c|l|c|}
\hline
Body $i$ &  mass  & \hspace{0.7cm}inertia tensor  & $\tau$ \\
\hline \hline
Body 1   &  78.10      & $\mathbb{I}^{(1)}=\textrm{diag}(0.08,26.05,26.1)$ & (0, 0, -1200)\\
\hline
Body 2   &  156.20     & $\mathbb{I}^{(2)}=\textrm{diag}(0.16,208.3,208.4)$ & -\\
\hline
Body 3   &  156.20     & $\mathbb{I}^{(3)}=\textrm{diag}(0.16,208.3,208.4)$ & -\\
\hline
\end{tabular}
\end{small}
\end{center}
\end{table}

\begin{table}[!htb]
\caption{Initial configuration of the planar quadrangle.}\label{TNNK:table2}
\begin{center}
\begin{small}
\begin{tabular}{|c||c|c|}
\hline
Body $i$ & $r^{(i)}$  & $\theta^{(i)}$\\
\hline
\hline
Body 1   & $\big(0.500, 0.866, 0\big)$ & $\big(0.866, 0, 0, 0.500\big)$\\
\hline
Body 2   & $\big(2.824, 2.553, 0\big)$ & $\big(0.978, 0, 0, 0.210\big)$\\
\hline
Body 3   & $\big(3.574, 1.687, 0\big)$ & $\big(0.877, 0, 0, 0.481\big)$\\
\hline
\end{tabular}
\end{small}
\end{center}
\end{table}

\begin{table}[!htb]
\caption{Joints of the planar quadrangle. Here $\chi^{(i)}$ is the position of the joint in $i^{th}$ body coordinate system.}
\label{TNNK:table3}
\begin{center}
\begin{small}
\begin{tabular}{|c||c|c|c|c|c|c|c|c|}
\hline
joint type &       $\mathcal{B}_i$         &    $\mathcal{B}_j$    & $\#$ of cond. & $\chi^{(i)}$ & $\chi^{(j)}$ &  $a^{(i)}$  &  $b^{(i)}$  &  $a^{(j)}$ \\
\hline
\hline
  revolute        &   0      &   1  &  5  & (0, 0, 0) & (-1, 0, 0) & $\mathsf{e}_2$ & $\mathsf{e}_1$ & $\mathsf{e}_3$\\
\hline
  spherical       &   1      &   2  &  3  & (1, 0, 0) & (-2, 0, 0) & -  & - & - \\
\hline
  cylindrical     &   2      &   3  &  4  & (2, 0, 0) & (2, 0, 0)  & $\mathsf{e}_1$ & $\mathsf{e}_2$ & $\mathsf{e}_3$\\
\hline
  revolute        &   3      &   0  &  5  & (-2, 0, 0) & $\left(2.5, 0, 0\right)$ & $\mathsf{e}_1$ & $\mathsf{e}_2$ & $\mathsf{e}_3$\\
\hline
\end{tabular}
\end{small}
\end{center}
\end{table}

Initial configuration is sketched in Figure \ref{fig:tasonivel_alkutila},
where global and local coordinate systems are also shown.
The local coordinate system of the ground body coincides with the global one.
The system starts from rest.

\subsubsection{Crank mechanism}

\begin{figure}[h!]
\begin{center}
\centering
\epsfig{file=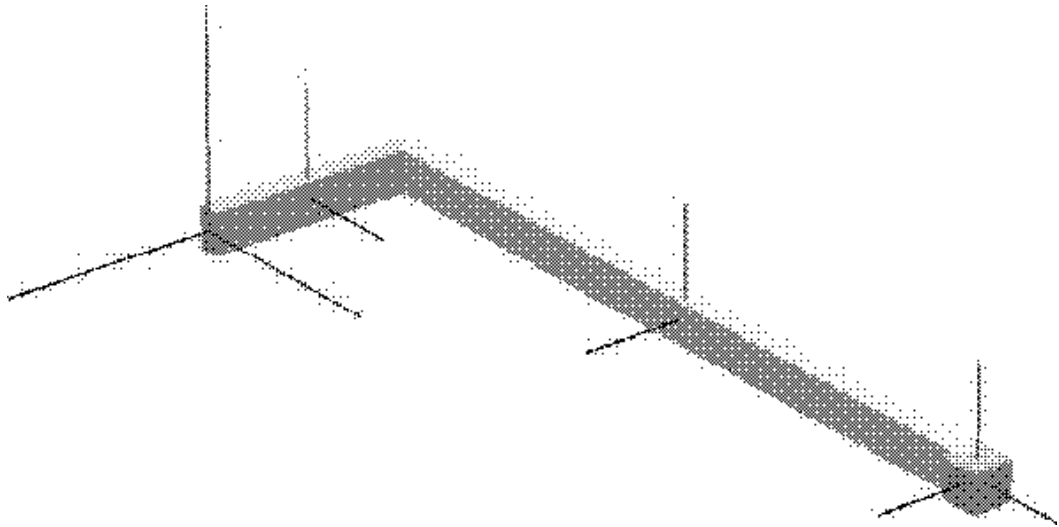,width=0.9\textwidth}
\caption{Initial configuration of the crank mechanism.}
\label{fig:kampimekanismi_alkutila}
\end{center}
\end{figure}

This example is of one degree of freedom, yet a real 3D-problem.
Again we have three solids, a ground body, and four joints: revolute,
spherical, universal, and translational.  There is an external torque
acting on the body 1 so that the corner joint is moving along a circle
centred at the global origin and perpendicular to $\mathsf{e}_1$.  The
spatial coordinate axes are in the Figure
\ref{fig:kampimekanismi_alkutila} as follows: $\mathsf{e}_1$ towards
the lower right corner, $\mathsf{e}_2$ upwards, $\mathsf{e}_3$ towards
the lower left corner.  Joint information is represented in Table
\ref{KM:table3}.  Parameters of the model are given in Tables
\ref{KM:table1} and \ref{KM:table2}. Initial configuration is
illustrated in Figure \ref{fig:kampimekanismi_alkutila} and the system
starts from rest.

\begin{table}[!htb]
\caption{Parameters of the crank mechanism.}
\label{KM:table1}
\begin{center}
\begin{small}
\begin{tabular}{|c||c|l|c|}
\hline
Body $i$ &  mass  &  \hspace{0.7cm}inertia tensor  & $\tau$ \\
\hline \hline
Body 1   &  19.50     & $\mathbb{I}^{(1)}=\textrm{diag}(0.02,0.41,0.42)$ & (0, 0, -50)\\
\hline
Body 2   &  70.29     & $\mathbb{I}^{(2)}=\textrm{diag}(0.07,18.99,19.04)$  & - \\
\hline
Body 3   &  7.81      & $\mathbb{I}^{(3)}=\textrm{diag}(0.01,0.01,0.01)$  & - \\
\hline
\end{tabular}
\end{small}
\end{center}
\end{table}

\begin{table}[!htb]
\caption{Initial configuration of the crank mechanism.}\label{KM:table2}
\begin{center}
\begin{small}
\begin{tabular}{|c||c|c|}
\hline
Body $i$ & $r^{(i)}$  & $\theta^{(i)}$ \\
\hline
\hline
Body 1   & $\big(0, 0, -0.25\big)$ & $\big(0.707, 0, 0.707, 0\big)$\\
\hline
Body 2   & $\big(0.9, 0, -0.5\big)$ & $\big(1, 0, 0, 0\big)$\\
\hline
Body 3   & $\big(1.8, 0 ,-0.5\big)$ & $\big(1, 0, 0, 0\big)$\\
\hline
\end{tabular}
\end{small}
\end{center}
\end{table}

\begin{table}[!htb]
\caption{Joints of the crank mechanism.}\label{KM:table3}
\begin{center}
\begin{small}
\begin{tabular}{|c||c|c|c|c|c|c|c|c|c|}
\hline
joint type &      $\mathcal{B}_i$     &    $\mathcal{B}_j$    & $\#$ of cond. &  $\chi^{(i)}$ & $\chi^{(j)}$ & $a^{(i)}$  &  $b^{(i)}$  &  $a^{(j)}$ & $b^{(j)}$\\
\hline
\hline
  spherical     &   1      &   2  &  3   &  $\left(0.25, 0, 0\right)$ &
$\left(-0.9, 0, 0\right)$   &  - & - & - & - \\
\hline
  translat.     &   0      &   3  &  5   &  (0, 0, 0)    &
$\left(-1.8, 0, 0.5\right)$ &  $\mathsf{e}_2$ & $\mathsf{e}_3$ & $\mathsf{e}_1$ & $\mathsf{e}_3$\\
\hline
  universal     &   2      &   3  &  4   &
$\left(0.9, 0, 0\right)$  & (0, 0, 0)    &
$\mathsf{e}_2$ & - & $\mathsf{e}_3$ & -\\
\hline
  revolute      &   1      &   0  &  5   &
$\left(-0.25, 0, 0\right)$ & (0, 0, 0)   &
$\mathsf{e}_2$ & $\mathsf{e}_1$ & $\mathsf{e}_1$ & -\\
\hline
\end{tabular}
\end{small}
\end{center}
\end{table}

\subsection{Test runs.}

Simulations were done with a 1.75GHz PC using \textsf{Linux} operating system, the algorithms were coded entirely in \textsf{C++} and compiled with the \textsf{Gnu C++} compiler. Figures \ref{fig:heiluri_alkutila},
\ref{fig:tasonivel_alkutila}, \ref{fig:kampimekanismi_alkutila} were
done with \textsf{OpenGL Framework} and other figures were plotted with
\textsf{Matlab}. In Tables \ref{PEND:results}, \ref{TNNK:results}, and \ref{KM:results}, CPU times are given in seconds.
Time interval was limited to $t=[0, 10]$ in all of the tests and the
absolute Newton tolerance was taken as  $\mathsf{tol_{abs}} = 10^{-5}$. 

We are mainly interested in the following comparisons:
\begin{itemize}
\item The quasi-orthogonal projection presented here (in the Algorithms \ref{algproj1} and  \ref{algprojholo}) vs. the usual orthogonal projection from \cite{arti14}.
\item The study of the effect of imposing/ignoring the energy conservation.
\end{itemize}
In the absence (resp. presence) of external forces we take $f_{\mathsf{inv}}=W$ (resp. $f_{\mathsf{inv}}=W_{\mathsf{total}}$). Recall that the computation of $W_{\mathsf{total}}$ involves a time integration which is handled with the trapezoidal rule.

\begin{figure}[h!]
\begin{center}
\centering
\includegraphics[width=8.5cm,height=6.5cm]{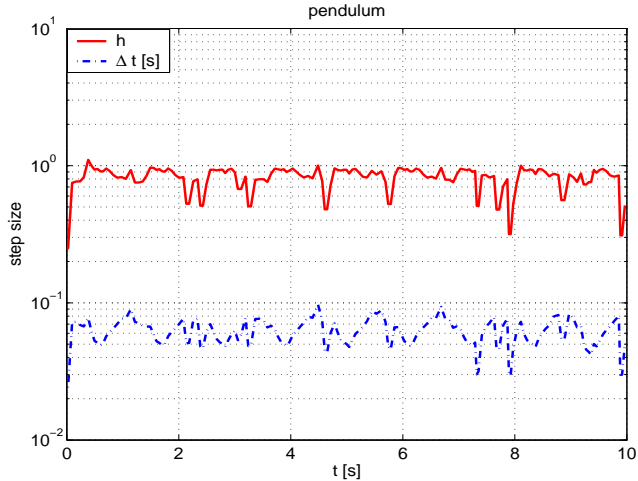}
\caption{Step sizes of the pendulum.}
\label{fig:PEND_askelpituudet}
\end{center}
\end{figure}

\noindent \textbf{3D pendulum.}  In the step size control we use
$\mathsf{Atol}=\mathsf{Rtol}=10^{-6}$, $\mathsf{fac_{max}} = 3.0$
\cite{hanowa}, and initial step size $h^0=0.25$.  The corresponding
step sizes are in Figure \ref{fig:PEND_askelpituudet} and run
characteristics in Table \ref{PEND:results}. 

Compared to our usual orthogonal projection \cite{arti14}, the
quasi-orthogonal projection needs 10\% more timesteps, and rejects
steps three times as much as the usual method, but uses a bit less
Newton iterations (on the average, yet the maximum number is twice
that of the usual orthogonal projection).  Still, the overall CPU time
spent by the new method is only a half, therefore the quasi-orthogonal
approach is 50\% faster.  This latter result illustrates well the
noticeable effectiveness of the quasi-orthogonal approach.  Then
again, the number of rejected steps and $\#df_{\mathsf{inv}}$,
$\#d^2f_{\mathsf{inv}}$ indicate some sort of sensitivity of the
quasi-orthogonal method, and it would be useful to find a strategy to
reduce this sensitivity making the method more robust in this sense.
Nevertheless, the results were qualitatively correct.

Also we like to point out that the quasi-orthogonal was more robust
than the usual orthogonal projection when an external torque $\tau$
was applied (to save space we do not tabulate these computations). 
When $|\tau|\approx 10$, including the energy invariant increased the computing time for both methods,
but the quasi-orthogonal one suffered much less ($110\%$ increase) than
the orthogonal one ($276\%$ increase).  When $|\tau|\approx 50$ the
quasi-orthogonal needed a lot more Newton iterations but the stepsize
stayed reasonable, whereas the orthogonal suffered from radically
decreasing stepsize.

In Table \ref{testit:ilman_energiaa}, by comparing the ``Pendulum''
panel to Table \ref{PEND:results}, can be seen that by ignoring the
energy invariant the computation has speeded up in both the
quasi-orthogonal ($20\%$) and the usual orthogonal projection ($5\%$)
as expected due to reduced amount of work.  But here more interesting
is that now the number of quasi-orthogonal projections is
significantly reduced: it is now about the same as (indeed even less
than) with the usual projection.  This shows that $f_{\mathsf{inv}}$
is, while itself a small system, slowing down the overall convergence of
the quasi-orthogonal method.  Nevertheless the quasi-orthogonal one is
58\% faster, the total CPU time is only a third of the usual method.

\begin{figure}[h!]
\begin{center}
\hspace{-1.5cm}
\begin{minipage}[t]{0.5\linewidth}
\centering
\includegraphics[width=8.5cm,height=6.5cm]{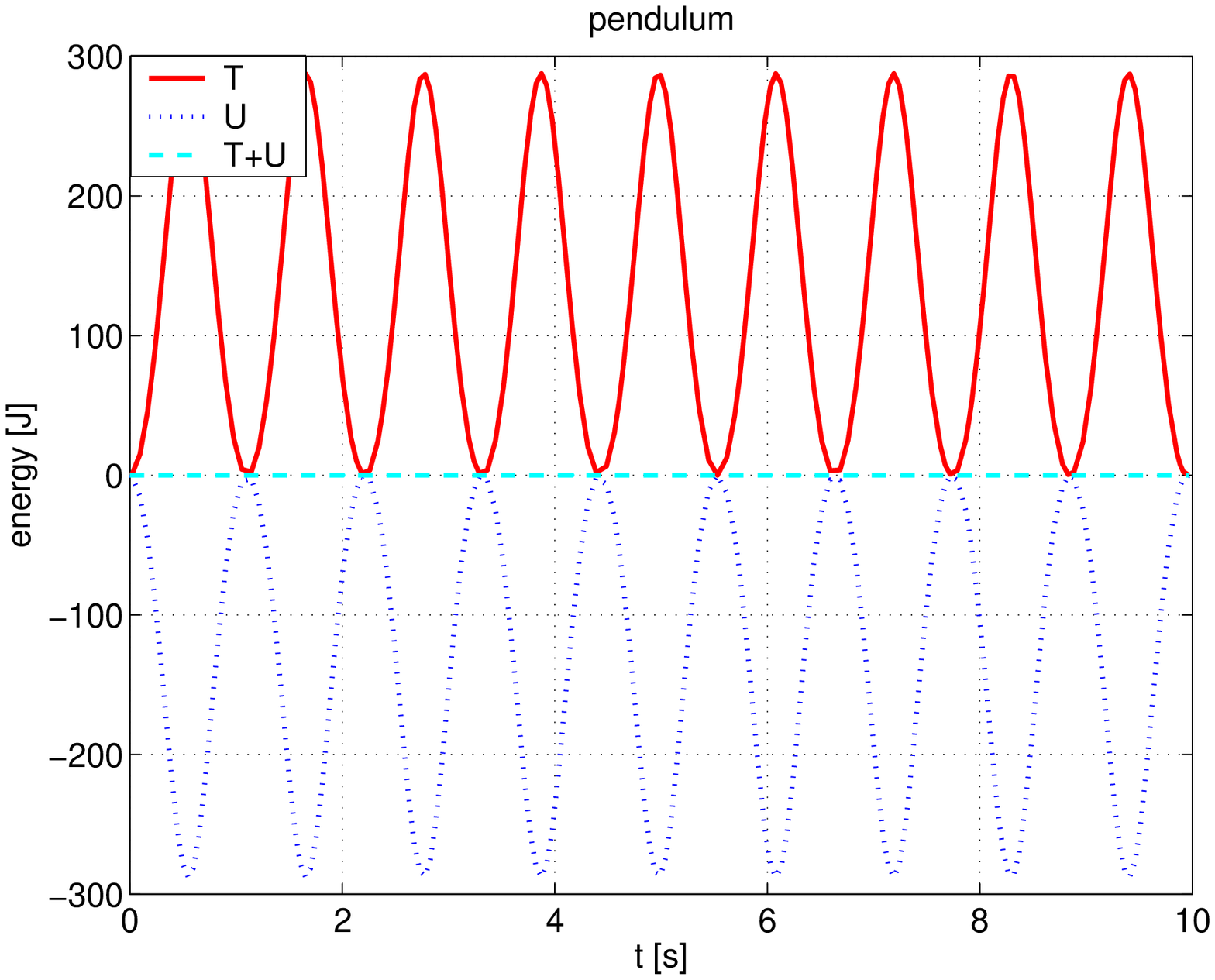}
\end{minipage}
\hspace{1.0cm}
\begin{minipage}[t]{0.5\linewidth}
\centering
\includegraphics[width=8.5cm,height=6.5cm]{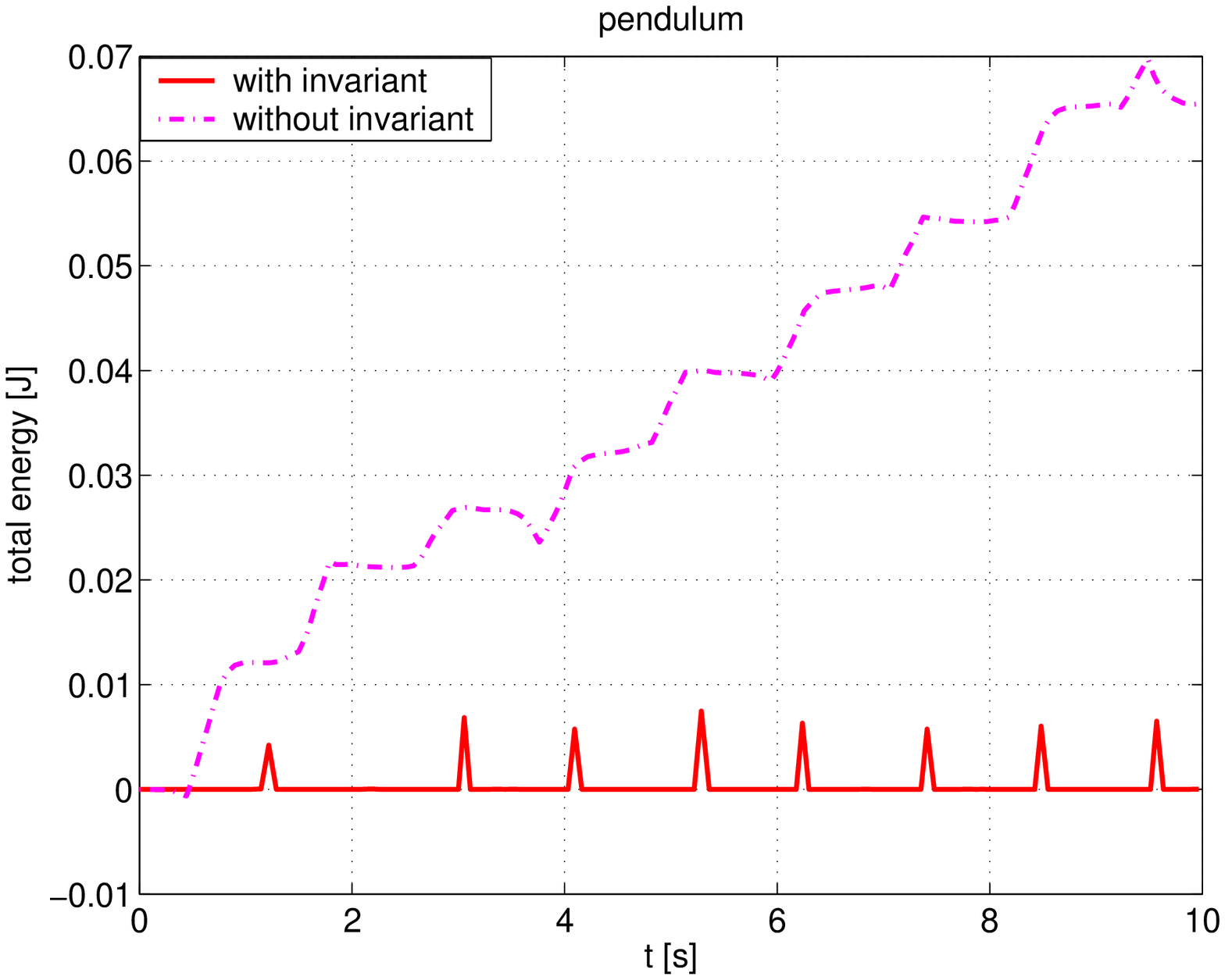}
\end{minipage}
\end{center}
\caption{On the left: energies of the pendulum when the conservation of the energy is imposed. On the right: the total energy with and without conservation of energy.}
\label{fig:PEND_energy}
\end{figure}

In Figure \ref{fig:PEND_energy} we have plotted fluctuations of different energies. Note that we do not need to do time integration to compute the total energy, and hence one source of numerical errors is eliminated. It is seen that the energy remains practically constant as it should. However, if the conservation of energy is not explictly imposed, the drift-off
is quite significant even on short time interval.

\begin{table}[!htb]
\caption{Run characteristics of the pendulum.}\label{PEND:results}
\begin{center}
\begin{tabular}{|c|c|c|}
\hline
  & quasi-orth.   & orthogonal \\
  & projection    & projection \\
\hline
$\textrm{$\#$ succ. (rej.) steps}$                   &   164(36)  & 151(12)  \\
\hline
Newton: av(max)                                      &   0.97(8) & 1.43(4)  \\
\hline
$\#$ dist                                            &   1395    & 1136  \\
$\#$ proj                                            &   1280    & 1070  \\
\hline
CPU dist                                             &   0.84    & 0.69  \\
CPU proj                                             &   1.65    & 4.30  \\
CPU total                                            &   2.53    & 5.05  \\
\hline
$\#~dg$               & 4856   &  6124 \\
$\#~d^2g^i$           & 9669  &  18372 \\
$\#~d^3g^i$           &   -    &  0 \\
$\#~df_{\mathsf{inv}}$       & 1211   &  631 \\
$\#~d^2f_{\mathsf{inv}}$     & 397   &  25 \\
\hline
\end{tabular}
\end{center}
\end{table}

\begin{figure}[h!]
\begin{center}
\hspace{-1.5cm}
\begin{minipage}[t]{0.5\linewidth}
\centering
\includegraphics[width=8.5cm,height=6.5cm]{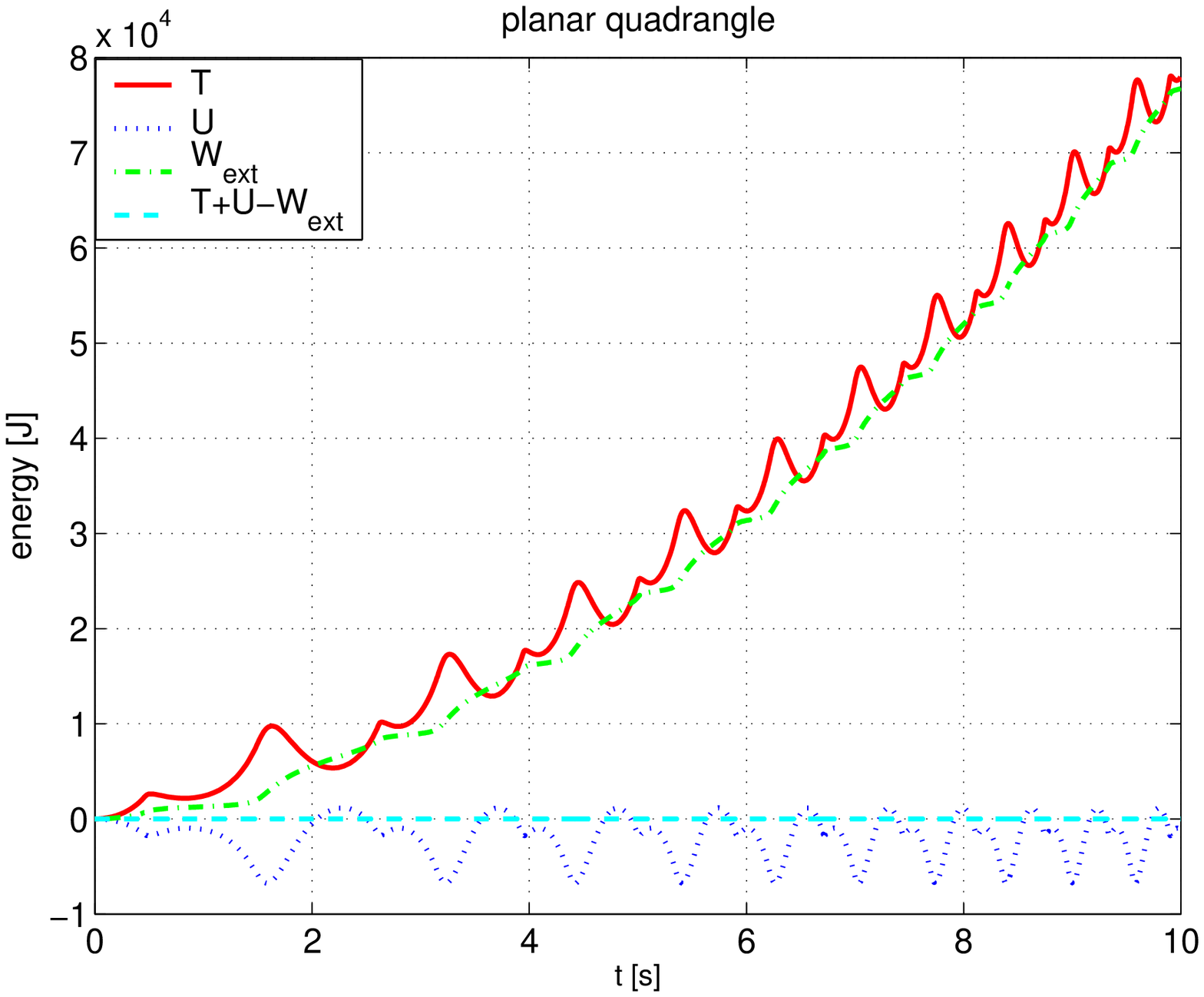}
\end{minipage}
\hspace{1.0cm}
\begin{minipage}[t]{0.5\linewidth}
\centering
\includegraphics[width=8.5cm,height=6.5cm]{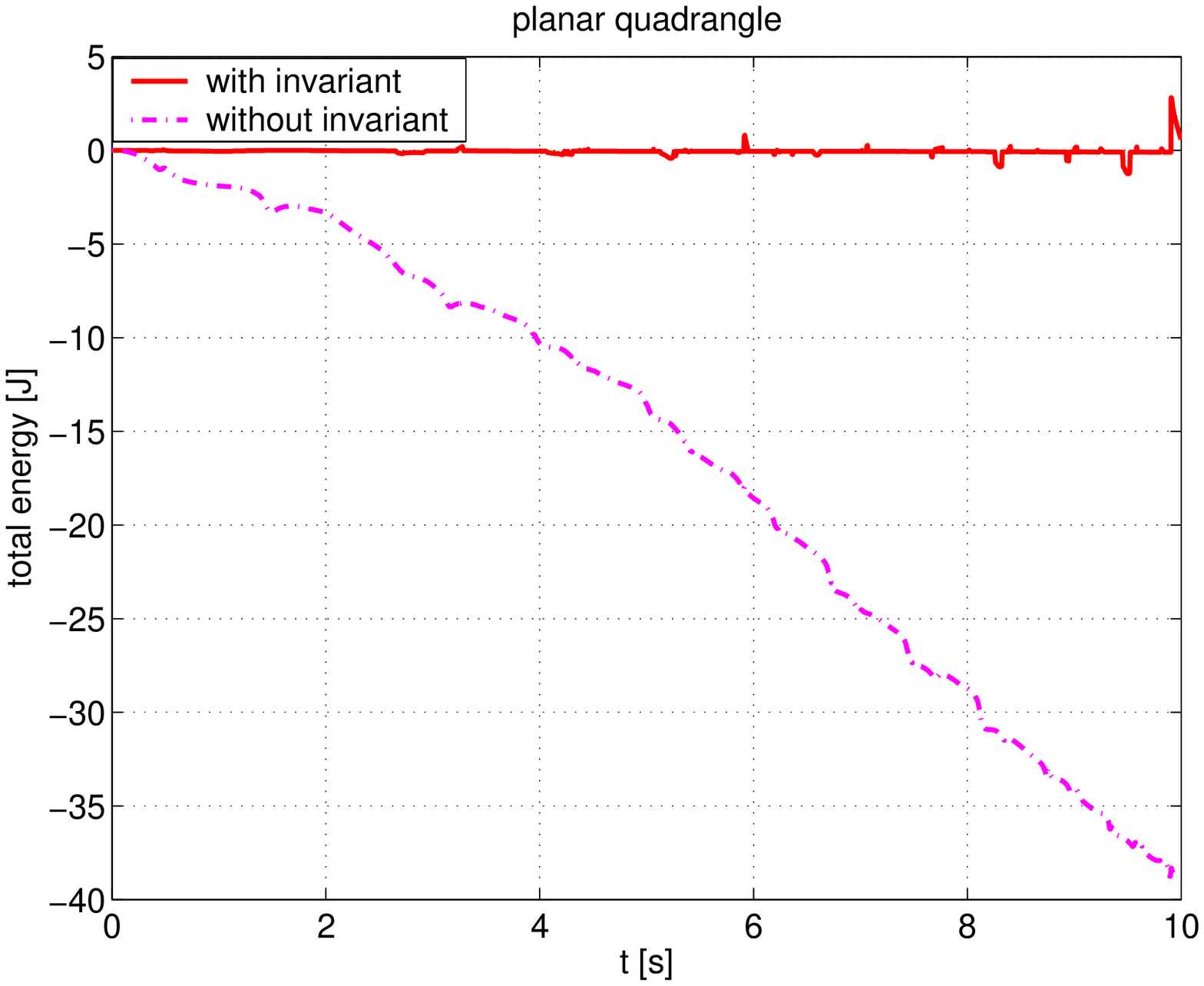}
\end{minipage}
\end{center}
\caption{On the left: energies of the planar quadrangle when the conservation of the energy is imposed. On the right: the total energy with and without conservation of energy.}
\label{fig:TNNK_energy}
\end{figure}

\begin{figure}[h!]
\begin{center}
\centering
\includegraphics[width=0.7\textwidth]{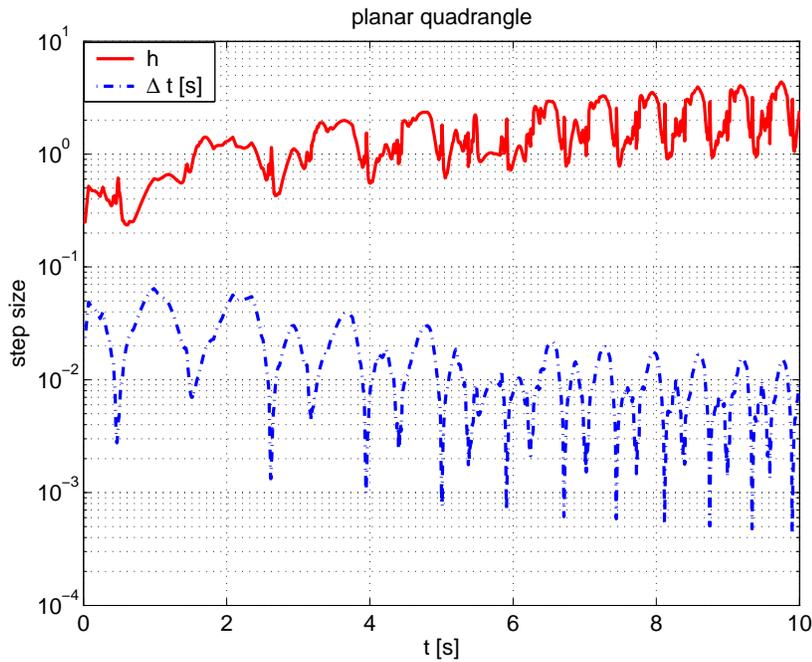}
\caption{Step sizes of the planar quadrangle.}
\label{fig:TNNK_askelpituudet}
\end{center}
\end{figure}

\noindent \textbf{Planar quadrangle.}  In the step size control we use
$\mathsf{Atol}=\mathsf{Rtol}=10^{-8}$, $\mathsf{fac_{max}} = 3.0$, and
initial step size $h^0=0.25$.  Run characteristics are in Table
\ref{TNNK:results}.  The step sizes and  the corresponding evolution of the energies with the quasi-orthogonal projection
are represented in Figures
\ref{fig:TNNK_askelpituudet} and \ref{fig:TNNK_energy}.
In the left panel of Figure \ref{fig:TNNK_energy} the potential energy is oscillating as one would physically expect since it consists
only of gravity, moreover the period of oscillation is getting shorter as the 
torque is speeding up the system.  This speedup is roughly 
linear in velocity and clearly visible in the 
quadratic tendency of the kinetic energy.

On the right panel is a magnification of the total energy which stays nearly constant when the conservation of energy is imposed. There are a few occasional spikes, that become more frequent and stronger with increasing kinetic energy.  The large number of maximum Newton iterations is related to these spikes. Without invariants the drift off of the energy is again quite rapid.

The quasi-orthogonal projection needs more steps (twice as many as the usual orthogonal projection) to converge, but is still $50\%$ faster in CPU time.  This shows the  quasi-orthogonal projection is 4 times faster.

\begin{table}[!htb]
\caption{Run characteristics of the planar quadrangle.}\label{TNNK:results}
\begin{center}
\begin{tabular}{|c|c|c|}
\hline
  & quasi-orth.   & orthogonal \\
  & projection    & projection \\
\hline
$\textrm{$\#$ succ. (rej.) steps}$                   &   1300(168)  & 1218(149)  \\
\hline
Newton: av(max)                                      &   3.25(265)  & 1.82(4)  \\
\hline
$\#$ dist                                            &   10139    &  9427 \\
$\#$ proj                                            &   9721     &  10932 \\
\hline
CPU dist                                             &   42.46    &  40.95 \\
CPU proj                                             &   175.66   &  457.73  \\
CPU total                                            &   218.91   &  499.41  \\
\hline
$\#~dg$               & 47842   &  57473 \\
$\#~d^2g^i$           & 613105  &  977041 \\
$\#~d^3g^i$           &   -     &  52388 \\
$\#~df_{\mathsf{inv}}$       & 27958   &  14961 \\
$\#~d^2f_{\mathsf{inv}}$     & 12194   &  4746 \\
\hline
\end{tabular}
\end{center}
\end{table}

\begin{figure}[h!]
\begin{center}
\hspace{-1.5cm}
\begin{minipage}[t]{0.5\linewidth}
\centering
\includegraphics[width=8.5cm,height=6.5cm]{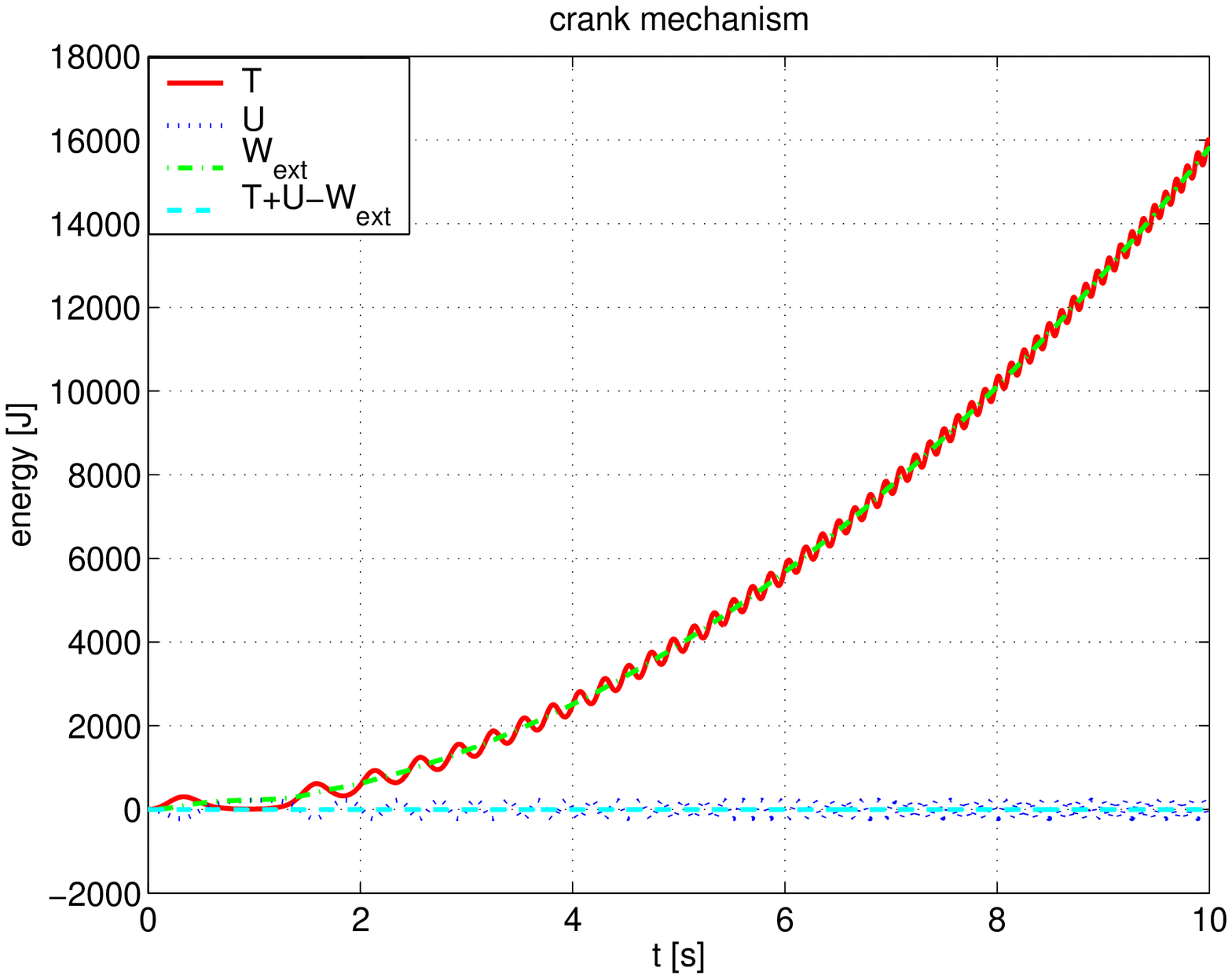}
\end{minipage}
\hspace{1.0cm}
\begin{minipage}[t]{0.5\linewidth}
\centering
\includegraphics[width=8.5cm,height=6.5cm]{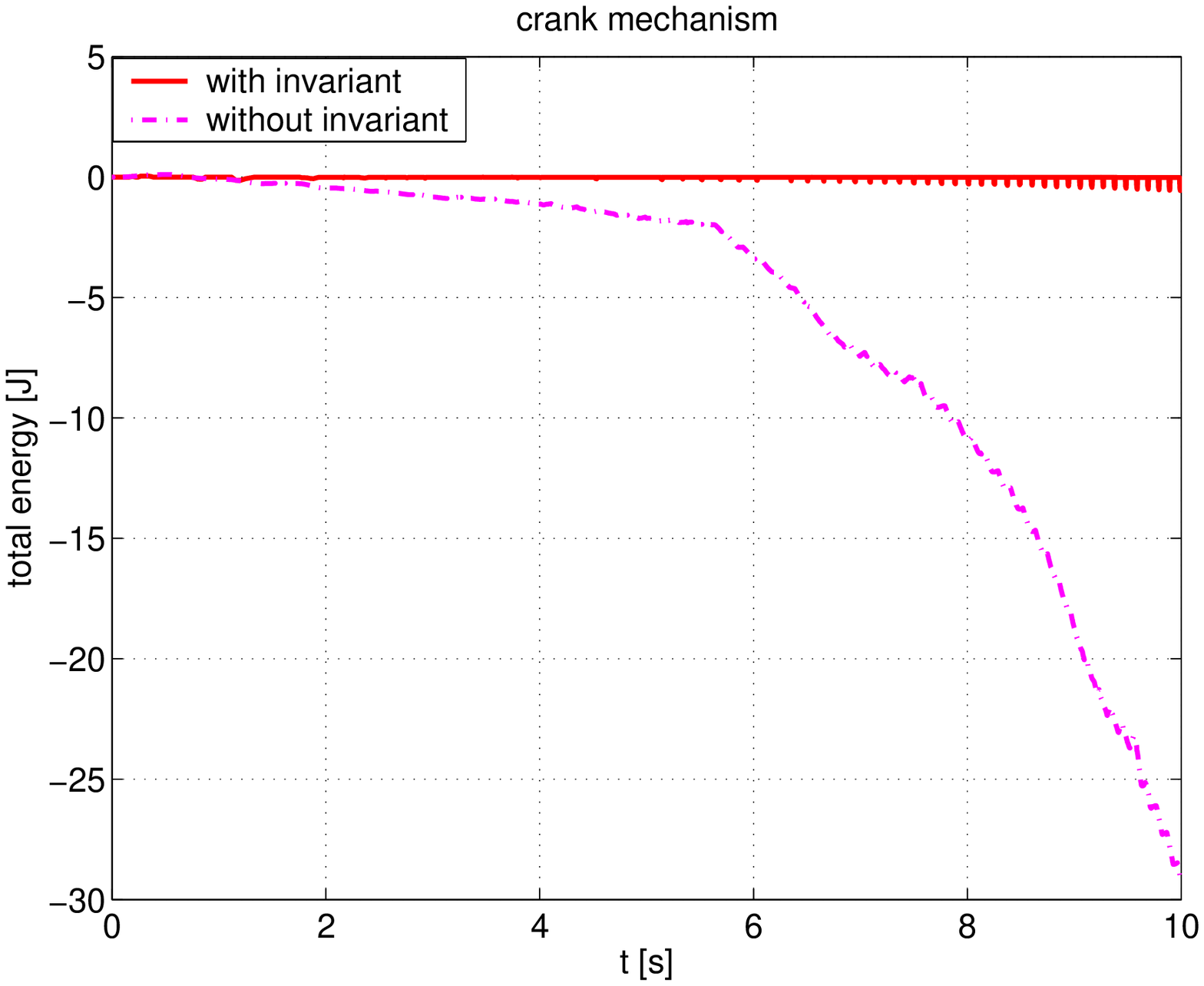}
\end{minipage}
\end{center}
\caption{On the left: energies of the crank mechanism when the conservation of the energy is imposed. On the right: the total energy with and without conservation of energy.}
\label{fig:KM_energy}
\end{figure}

\begin{figure}[h!]
\begin{center}
\centering
\includegraphics[width=8.5cm,height=6.5cm]{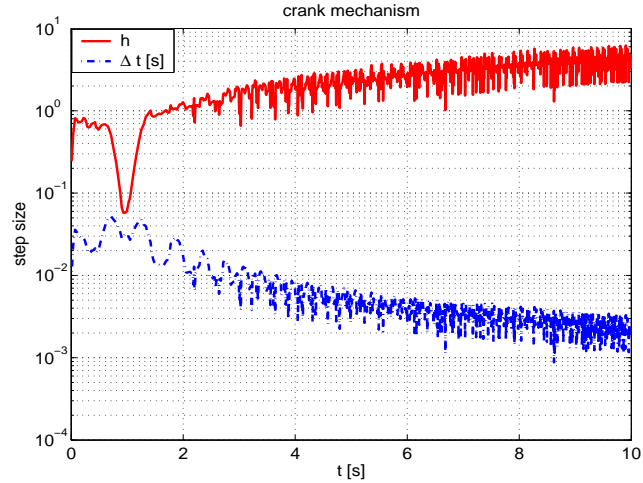}
\caption{Step sizes of the crank mechanism.}
\label{fig:KM_askelpituudet}
\end{center}
\end{figure}

\noindent \textbf{Crank mechanism.}  In the step size control we use
$\mathsf{Atol}=\mathsf{Rtol}=10^{-7}$, $\mathsf{fac_{max}} = 3.0$, and
initial step size $h^0=0.25$.  Run characteristics are in Table
\ref{KM:results}.  Evolution of energies is represented in
Figure \ref{fig:KM_energy} and the step sizes in Figure
\ref{fig:KM_askelpituudet}.  The results are similar to the planar
quadrangle case so we will be brief here.  The most notable
differences are that the differentials of $f_{\mathsf{inv}}$ are evaluated 3-5
times as often, and the quasi-orthogonal projection needs about the
same number of steps yet uses only a quarter of CPU time compared to
the usual orthogonal projection.  Hence the quasi-orthogonal
projection is 4 times faster here as well.

\begin{table}[!htb]
\caption{Run characteristics of the crank mechanism.}\label{KM:results}
\begin{center}
\begin{tabular}{|c|c|c|}
\hline
  & quasi-orth.   & orthogonal \\
  & projection    & projection \\
\hline
$\textrm{$\#$ succ. (rej.) steps}$                   &   2086(352)  &  2070(113) \\
\hline
Newton: av(max)                                      &   6.69(215)  &  2.61(6) \\
\hline
$\#$ dist                                            &   16859    &  15175 \\
$\#$ proj                                            &   13358    &  17458 \\
\hline
CPU dist                                             &   88.08    &  78.44 \\
CPU proj                                             &   757.65   &  3198.90  \\
CPU total                                            &   848.67   &  3292.00  \\
\hline
$\#~dg$               & 124303  & 123597  \\
$\#~d^2g^i$           & 1778438 & 2101149  \\
$\#~d^3g^i$           &   -     & 36758  \\
$\#~df_{\mathsf{inv}}$       & 108843  & 31219  \\
$\#~d^2f_{\mathsf{inv}}$     & 58003   & 10981  \\
\hline
\end{tabular}
\end{center}
\end{table}

\begin{table}[!htb]
\caption{Test runs without the conservation of energy.}\label{testit:ilman_energiaa}
\begin{center}
\begin{tabular}{|l|c|c|}
\hline
  & quasi-orth.   & orthogonal \\
  & projection    & projection \\
\hline
pendulum: &  &  \\
\hspace{1cm}$\textrm{$\#$ succ. (rej.) steps}$             &  153(23)   &  152(13)  \\
\hspace{1cm}Newton: av(max)                                &  0.929(2)  &  0.987(4) \\
\hspace{1cm}$\#$ proj                                      &  795       &  833  \\
\hspace{1cm}CPU total                                      &  1.98      &  4.80 \\
\hspace{1cm}$\#~dg$                                        &  4158      &  5885 \\
\hline
planar quadrangle: &  &  \\
\hspace{1cm}$\textrm{$\#$ succ. (rej.) steps}$             &  1487(181) &  1226(161) \\
\hspace{1cm}Newton: av(max)                                &  0.113(2)  &  0.995(3)  \\
\hspace{1cm}$\#$ proj                                      &  4332      &  6852  \\
\hspace{1cm}CPU total                                      &  95.03     &  301.05 \\
\hspace{1cm}$\#~dg$                                        &  32003     &  40476 \\
\hline
crank mechanism: &  &  \\
\hspace{1cm}$\textrm{$\#$ succ. (rej.) steps}$             &  1996(87)  &  2077(129) \\
\hspace{1cm}Newton: av(max)                                &  0.325(2)  &  1.42(5)   \\
\hspace{1cm}$\#$ proj                                      &  7133      &  11147  \\
\hspace{1cm}CPU total                                      &  191.27    &  2099.30  \\
\hspace{1cm}$\#~dg$                                        &  45054     &  87018 \\
\hline
\end{tabular}
\end{center}
\end{table}

\section{Conclusion and perspectives}
\label{sec:concl}
We have derived a computational model to simulate multibody dynamics with holonomic constraints.  This approach is based on jet spaces and Lagrangian formalism and avoids the traditional drift-off problems by an orthogonal projection onto the relevant manifold.

As we have seen in the numerical examples above our method can take into account arbitrary (holonomic) constraints and in addition we can include any invariants, such as the conservation of energy, in our model. Hence our simulations can run indefinitely as far as the physical relevance of the constraints is concerned.

Computationally the most expensive part in our simulations is the projection step; this may take as much as 90\% of the total time. However, this aspect is not fully optimized in our code. Much of the time goes to updating various differentials, and it is clear that these updates could be less frequent. Another possibility to speed up the code would be to use \emph{automatic differentiation} \cite{griewank}. However, exploring this idea was left to future work.

Another topic that we aim to work on is to improve further the
quasi-orthogonal iteration: now we iterate first onto $M_{\mathsf{hc}}
\cap M_{\mathsf{inv}}$ and then onto $M_{\theta}$.  However, during
the latter step $\theta_1$ may change significantly and we
need to re-iterate onto $M_{\mathsf{hc}} \cap M_{\mathsf{inv}}$ again.
One possible reason for this is that in some cases the condition number of the relevant matrix in the Newton iteration was relatively big, resulting in the slow convergence. This could probably be fixed by some suitable precondition method.
In any case it would be useful to develop a strategy to get onto the
$M_{\theta}\cap M_{\mathsf{hc}} \cap M_{\mathsf{inv}}$ without losing
the efficiency of the quasi-orthogonal iteration.

We did not consider nonholonomic systems in the present article because treating them would have augmented the length of the paper considerably. The formulation of nonholonomic problems is very different from the holonomic ones, for example nonholonomic problems are \emph{not} variational (in a standard sense); see \cite[p. 208]{bloch} for a discussion and further references. However, we believe that the general framework of jet spaces is also suitable for numerical solution of nonholonomic systems and we hope to treat this case in future papers.

\bibliographystyle{amsplain}

\begin{thebibliography}{10}

\bibitem{amirouche}
F.~Amirouche, \emph{Fundamentals of multibody dynamics}, Birkh\"auser Boston
  Inc., Boston, MA, 2006.

\bibitem{appell}
P.~Appell, \emph{Trait\'e de m\'ecanique rationelle, \emph{ {T}ome {I} \&
  {T}ome {II}}}, \'Editions Jacques Gabay, 1991, r\'eimpression de $6^e$ \'ed.
  publi\'ee par Gauthier-Villars en 1941 ({T}ome {I}) et en 1953 ({T}ome {II}).

\bibitem{arnold}
V.~I. Arnold, \emph{Mathematical methods of classical mechanics}, 2nd ed.,
  Graduate Texts in Mathematics, vol.~60, Springer-Verlag, New York, 1989.

\bibitem{baumgarte}
J.~Baumgarte, \emph{Stabilization of constraints and integrals of motion in
  dynamical systems}, Comput. Methods Appl. Mech. Engrg. \textbf{1} (1972),
  1--16.

\bibitem{bloch}
A.~Bloch, \emph{Nonholonomic mechanics and control}, Interdisciplinary Applied
  Mathematics, vol.~24, Springer-Verlag, New York, 2003.

\bibitem{bullew}
F.~Bullo and A.~Lewis, \emph{Geometric control of mechanical systems}, Texts in
  applied mathematics, vol.~49, Springer, 2005.

\bibitem{gdjbayo}
J.~Garc{\'{\i}}a de~Jal{\'o}n and E.~Bayo, \emph{Kinematic and dynamic
  simulation of multibody systems}, Mechanical Engineering Series,
  Springer-Verlag, New York, 1994.

\bibitem{deeist}
R.~Dembo, S.~Eisenstat, and T.~Steinhaug, \emph{Inexact {N}ewton methods}, SIAM
  J. Numer. Anal. \textbf{19} (1982), no.~2, 400--408.

\bibitem{giahil}
M.~Giaquinta and S.~Hildebrandt, \emph{Calculus of variations {I}},
  Grundlehren, vol. 310, Springer, 1996.

\bibitem{griewank}
A.~Griewank, \emph{Evaluating derivatives}, SIAM, 2000.

\bibitem{hanowa}
E.~Hairer, S.~N\o rsett, and G.~Wanner, \emph{Solving ordinary differential
  equations {I}, {N}onstiff problems}, 2nd ed., Springer series in comp. math.,
  vol.~8, Springer, 1993.

\bibitem{krupkova}
O.~Krupkov{\'a}, \emph{Mechanical systems with nonholonomic constraints}, J.
  Math. Phys. \textbf{38} (1997), no.~10, 5098--5126.

\bibitem{lanczos}
C.~Lanczos, \emph{The variational principles of mechanics}, 4th ed., Dover,
  1970.

\bibitem{pommaret}
J.F. Pommaret, \emph{Systems of partial differential equations and {L}ie
  pseudo\-groups}, Gordon \&{} Breach, London, 1978.

\bibitem{rabrhe}
P.~Rabier and W.~Rheinboldt, \emph{Nonholonomic motion of rigid mechanical
  systems from a {DAE} viewpoint}, SIAM, 2000.

\bibitem{saunders}
D.~Saunders, \emph{The geometry of jet bundles}, London Math. Soc. Lecture note
  series, vol. 142, Cambridge university press, 1989.

\bibitem{schwerin}
R.~von Schwerin, \emph{Multibody system simulation}, Lecture Notes in
  Computational Science and Engineering, vol.~7, Springer-Verlag, Berlin, 1999.

\bibitem{seiler}
W.M. Seiler, \emph{Involution --- the formal theory of differential equations
  and its applications in computer algebra and numerical analysis},
  Habilitation thesis, Dept.\ of Mathematics, Universit\"at Mannheim, 2001,
  (manuscript accepted for publication by Springer-Verlag).

\bibitem{shabana}
A.~A. Shabana, \emph{Dynamics of multibody systems}, 2nd ed., Cambridge
  University Press, 1998.

\bibitem{spivak}
M.~Spivak, \emph{A comprehensive introduction to differential geometry, vol 1 -
  5}, 2nd ed., Publish or Perish, 1979.

\bibitem{tarkhanov}
N.~N. Tarkhanov, \emph{Complexes of differential operators}, Mathematics and
  its Applications, vol. 340, Kluwer Academic Publishers Group, Dordrecht,
  1995.

\bibitem{julk31}
J.~Tuomela and T.~Arponen, \emph{On the numerical solution of involutive
  ordinary differential systems}, IMA J. Num. Anal. \textbf{20} (2000),
  561--599.

\bibitem{julk36}
\bysame, \emph{On the numerical solution of involutive ordinary differential
  systems: {H}igher order methods}, BIT \textbf{41} (2001), 599--628.

\bibitem{arti14}
J.~Tuomela, T.~Arponen, and V.~Normi, \emph{On the numerical solution of
  involutive ordinary differential systems: {E}nhanced linear algebra}, to
  appear in IMA J. Num. Anal.

\end{thebibliography}
\providecommand{\bysame}{\leavevmode\hbox to3em{\hrulefill}\thinspace}
\providecommand{\MR}{\relax\ifhmode\unskip\space\fi MR }
\providecommand{\MRhref}[2]{%
  \href{http://www.ams.org/mathscinet-getitem?mr=#1}{#2}
}
\providecommand{\href}[2]{#2}

\end{document}